\documentclass[11pt, final, journal, onecolumn]{IEEEtran} 

\usepackage{times}
\usepackage{hyperref}
\usepackage{url}
\usepackage{mathtools} 
\usepackage{cases}
\usepackage{subcaption}
\usepackage{amssymb,algorithm,cite,caption,cases,float,graphicx,url,color}
\usepackage{thm-restate, enumerate }

\definecolor{darkred}{RGB}{250,0,0}
\definecolor{darkgreen}{RGB}{0,150,0}
\definecolor{myblue}{RGB}{0,0,250}
\definecolor{darkblue}{RGB}{0,0,200}
\hypersetup{colorlinks=true, linkcolor=darkred, citecolor=myblue, urlcolor=darkblue}

\usepackage{microtype} 

\usepackage[hmarginratio=1:1,top=32mm,columnsep=20pt]{geometry} 
\newtheorem{defn}{Definition}[section]

\newtheorem{thm}{Theorem}[section]

\newtheorem{lem}{Lemma}[section]

\newtheorem{ass}{Assumption}

\newcommand{\Roc}{\cite[Cor.~37.3.2]{Roc70} }
\newcommand{\consist}{\cite[Thm.~2.7]{NF36} }

\newcommand{\f}{{f}}
\newcommand{\Pp}{\mathbf{P}^\bot}
\newcommand{\al}{\alpha}

\newcommand{\lam}{\la_{\text{min}}}
\newcommand{\Rcrit}{\mathcal{R}_{\text{crit}}}
\newcommand{\Rbad}{\mathcal{R}_{\text{bad}}}
\newcommand{\ellbLM}{\ellb^{LM}}
\newcommand{\tbLM}{\tb^{LM}}
\newcommand{\ah}{\hat\alpha}
\newcommand{\vbs}{\vb_*}
\newcommand{\eps}{\epsilon}
\newcommand{\hb}{\bar\h}
\newcommand{\fb}{\bar{f}^*}
\newcommand{\wt}{\tilde\w}
\newcommand{\lin}{\text{lin}}
\newcommand{\sign}{\mathrm{sign}}
\newcommand{\DK}{D_{\Kc}(\mu\x_0)}
\newcommand{\etav}{\vec{\eta}}

\newcommand{\kac}{\kappa_\text{crit}}

\newcommand{\prox}{\mathrm{prox}}

\newcommand{\tauu}{\tau^2}
\newcommand{\one}{\mathbf{1}}
\newcommand{\I}{\mathbf{I}}
\newcommand{\gv}{\vec{g}}
\newcommand{\xx}{\overline{\x}}

\newcommand{\XX}{\overline{X}}

\newcommand{\mn}{{(n)}}
\newcommand{\E}{\mathbb{E}}                    
\newcommand{\la}{{\lambda}}                     
\newcommand{\sigg}{\sigma^2}
\newcommand{\sig}{\sigma}

\newcommand{\vecc}{\mathrm{vec}}

\newcommand{\nn}{\notag}
\newcommand{\R}{\mathbb{R}}

\newcommand{\Ub}{\mathbf{U}}

\newcommand{\Hb}{\mathbf{H}}
\newcommand{\G}{\mathbf{G}}
\newcommand{\Sb}{\mathbf{S}}
\newcommand{\X}{\mathbf{X}}
\newcommand{\XXb}{\overline{\mathbf{X}}}
\newcommand{\A}{\mathbf{A}}

\newcommand{\Vb}{\mathbf{V}}

\newcommand{\ellb}{\boldsymbol{\ell}}
\newcommand{\x}{\mathbf{x}}
\newcommand{\w}{\mathbf{w}}
\newcommand{\ub}{\mathbf{u}}
\newcommand{\g}{\mathbf{g}}
\newcommand{\vb}{\mathbf{v}}

\newcommand{\e}{\mathbf{e}}
\newcommand{\tb}{\mathbf{t}}
\newcommand{\y}{\mathbf{y}}
\newcommand{\s}{\mathbf{s}}
\newcommand{\z}{\mathbf{z}}

\newcommand{\ab}{\mathbf{a}}

\newcommand{\h}{\mathbf{h}}

\newcommand{\Sc}{{\mathcal{S}}}

\newcommand{\Kc}{\mathcal{K}}

\newcommand{\Nn}{\mathcal{N}}

\newcommand{\Qc}{\mathcal{Q}}

\newcommand{\beq}{\begin{equation}}
\newcommand{\eeq}{\end{equation}}
\newcommand{\bea}{\begin{align}}
\newcommand{\eea}{\end{align}}

\newcommand{\lac}{\lambda_{\text{crit}}}

\newcommand{\order}[1]{\mathcal{O}\left(#1\right)}

\newcommand{\rP}{\xrightarrow{P}}

\newcommand{\dom}{\operatorname{dom}}

\title{
The LASSO with Non-linear Measurements is Equivalent to One With Linear Measurements
}

\author{\vspace{10pt}
Christos Thrampoulidis,\hspace{10pt}Ehsan Abbasi,\hspace{10pt}Babak Hassibi\vspace{10pt}
\thanks{Department of Electrical Engineering, Caltech, Pasadena -- 91125, 
emails: {\tt(cthrampo, eabbasi, hassibi)@caltech.edu}. }
}
\begin{document}

\maketitle

\begin{abstract}
Consider estimating an unknown, but structured (e.g. sparse, low-rank, etc.), signal $\mathbf{x}
_0\in\mathbb{R}^n$ from a vector $\mathbf{y}\in\mathbb{R}^m$ of measurements of the form $y_i=g_i({\mathbf{a}_i}^T\mathbf{x}_0)$, where the ${\mathbf{a}_i}$'s are the rows of a known measurement matrix $\mathbf{A}$, and, $g(\cdot)$ is a (potentially unknown) nonlinear and random link-function. Such measurement functions could arise in applications where the measurement device has nonlinearities and uncertainties. It could also arise by design, e.g., $g_i(x)=\mathrm{sign}(x+z_i)$, corresponds to noisy 1-bit quantized measurements. Motivated by the classical work of Brillinger, and more recent work of Plan and Vershynin, we estimate $\mathbf{x}_0$ via solving the Generalized-LASSO, i.e., $\hat{\mathbf{x}}:=\arg\min_{\mathbf{x}}\|\mathbf{y}-\mathbf{A}\mathbf{x}_0\|_2+\lambda f(\mathbf{x})$ for some regularization parameter $\lambda>0$ and some (typically non-smooth) convex regularizer $f(\cdot)$ that promotes the structure of $\mathbf{x}_0$, e.g. $\ell_1$-norm, nuclear-norm, etc. While this approach seems to naively ignore the nonlinear function $g(\cdot)$, both Brillinger (in the non-constrained case) and Plan and Vershynin have shown that, when the entries of $\mathbf{A}$ are iid standard normal, this is a good estimator of $\mathbf{x}_0$ up to a constant of proportionality $\mu$, which only depends on $g(\cdot)$. In this work, we considerably strengthen these results by obtaining explicit expressions for$\|\hat{\mathbf{x}}-\mu\mathbf{x}_0\|_2$, for the \emph{regularized} Generalized-LASSO, that are asymptotically \emph{precise} when $m$ and $n$ grow large. A main result is that the estimation performance of the Generalized LASSO with non-linear measurements is \emph{asymptotically the same} as one whose measurements are linear $y_i=\mu {\mathbf{a}_i}^T\mathbf{x}_0 + \sigma z_i$, with $\mu = \mathbb{E}\gamma g(\gamma)$ and $\sigma^2 = \mathbb{E}(g(\gamma)-\mu\gamma)^2$, and, $\gamma$ standard normal. To the best of our knowledge, the derived expressions on the estimation performance are the first-known precise results in this context. One interesting consequence of our result is that the optimal quantizer of the measurements that minimizes the estimation error of the Generalized LASSO is the celebrated Lloyd-Max quantizer.

\end{abstract}
%

\section{Introduction}
\subsection{Problem Setup}

\subsubsection{Non-linear Measurements}
Consider the problem of estimating an unknown signal vector $\x_0\in\R^n$ from a vector $\y=(y_1,y_2,\ldots,y_m)^T$ of $m$ measurements taking the following form:
\begin{align}\label{eq:model}
y_i = g_i(\ab_i^T\x_0), \quad i=1,2,\ldots,m.
\end{align}
Here, each $\ab_i$ represents a (known) measurement vector. 
The $g_i$'s are independent copies of a (generically random) link function $g$. For instance, $g_i(x)=x+z_i$, with say $z_i$ being normally distributed, recovers the standard linear regression setup with gaussian noise. In this paper, we are particularly interested in scenarios where $g$ is \emph{non-linear}. Notable examples include $g(x) = \sign(x)$ (or $g_i(x)=\sign(x+z_i)$) and $g(x)=(x)_+$, corresponding to $1$-bit quantized (noisy) measurements, and, to the censored Tobit model, respectively. Depending on the situation, $g$ might be known or unspecified. In the statistics and econometrics literature, the measurement model in \eqref{eq:model} is popular under the name \emph{single-index model} and several aspects of it have been well-studied, e.g.\cite{Bri,Bri2,ichimura1993semiparametric,li1989regression}. 

%

%
\vspace{8pt}
\subsubsection{Structured Signals}
It is typical in many instances that the unknown signal $\x_0$ obeys some sort of \emph{structure}. For instance, it might be sparse in which case only a few $k\ll n$, of its entries are non-zero; or, it might be that $\x_0=\vecc(\X_0)$, where $\X_0\in\R^{\sqrt{n}\times\sqrt{n}}$ is a matrix of low-rank $r\ll n$. 
 To exploit this information it is typical to associate with the structure of $\x_0$ a properly chosen function $f:\R^n\rightarrow\R$, which we  refer to as the \emph{regularizer}. Of particular interest are \emph{convex} and non-smooth such regularizers, e.g. the $\ell_1$-norm for sparse signals, the nuclear-norm for low-rank ones, etc. Please refer for example to \cite{Cha,bach2010structured,halabi2014totally,TroppEdge} for further discussions.
 
 \vspace{8pt}
 \subsubsection{An Algorithm for Linear Measurements: The Generalized LASSO}
 When the link function is \emph{linear}, i.e. $g_i(x)=x+z_i$,
perhaps the most popular way of estimating $\x_0$ is via solving the Generalized LASSO 
 algorithm:
 \begin{align}\label{eq:LASSO1}
 \hat\x := \arg\min_{\x}\|\y-\A\x\|_2 + \la f(\x).
 \end{align}
 Here, $\A=[\ab_1, \ab_2, \ldots, \ab_m]^T\in\R^{m\times n}$ is the known measurement matrix and $\la>0$ is a regularizer parameter. This is often referred to as the $\ell_2$-LASSO or the square-root-LASSO \cite{Belloni} to distinguish from the one which solves $\min_{\x} \frac{1}{2}\|\y-\A\x\|_2^2 + \la f(\x)$, instead. The results of this paper can be accustomed to this latter version, but for concreteness, we restrict  attention to \eqref{eq:LASSO1} throughout. The acronym LASSO for \eqref{eq:LASSO1} was introduced 
 in \cite{TibLASSO} for the special case of $\ell_1$-regularization; \eqref{eq:LASSO1} is a natural generalization to other kinds of structures and includes the group-LASSO\cite{groupLASSO}, the fused-LASSO\cite{fusedLASSO} as special cases. We often drop the term ``Generalized" and refer to \eqref{eq:LASSO1} simply as the LASSO.  
 
 One popular,  measure of  estimation performance of \eqref{eq:LASSO1} is the squared-error $\|\hat\x-\x_0\|_2^2$.  Recently, there have been significant advances on establishing tight bounds and even \emph{precise} characterizations of this quantity, in the presence of linear measurements \cite{DMM,montanariLASSO,StoLASSO,OTH13,ISIT15_ell22,COLT15}. 
 Such precise results have been core to building a better understanding of the behavior of the LASSO, and, in particular, on the exact role played by the choice of the regularizer $f$ (in accordance with the structure of $\x_0$), by the number of measurements $m$, by the value of $\la$, etc.. In certain cases, they even provide us with useful insights into practical matters such as the tuning of the regularizer parameter. 
 
 
 \vspace{5pt}
\vspace{8pt}
 \subsubsection{Using the LASSO for Non-linear Measurements?}
 The  LASSO is by nature tailored to a linear model for the measurements. Indeed, the first term of the objective function in \eqref{eq:LASSO1} tries to fit  $\A\x$ to the observed vector $\y$ presuming that this is of the form $y_i=\ab_i^T\x_0+ \text{noise}$. Of course, no one stops us from continuing to use it even in cases where $y_i=g(\ab_i^T\x_0)$ with $g$ being \emph{non}-linear\footnote{Note that
  the Generalized LASSO in \eqref{eq:LASSO1} does not assume knowledge of $g$. All that is assumed is the availability of the measurements $y_i$. Thus, the link-function might as well be unknown or unspecified.}. But, the question then becomes: 
 Can there be any guarantees that the solution $\hat\x$ of the Generalized LASSO is still a good estimate of $\x_0$? 



The question just posed was first studied back in the early 80's by Brillinger \cite{Bri} who provided answers in the case of solving \eqref{eq:LASSO1} without a regularizer term. This, of course, corresponds to  standard Least Squares (LS). Interestingly, he showed that when the measurement vectors are Gaussian, then the LS solution is a consistent estimate of $\x_0$, up to a constant of proportionality $\mu$, which only depends on the link-function $g$. The result is sharp, but only  under the assumption that the number of measurements $m$ grows large, while the signal dimension $n$ stays fixed, which was the typical setting of interest at the time. 
In the world of structured signals and high-dimensional measurements, the problem was only 
very recently revisited by Plan and Vershynin \cite{Ver}. They consider a \emph{constrained} version of the Generalized LASSO, in which the regularizer is essentially replaced by a constraint, and derive upper bounds on its performance. The bounds are not tight (they involve absolute constants), but they demonstrate some key features: i) the solution to the constrained LASSO $\hat\x$ is a good estimate of $\x_0$ up to the  same constant of proportionality $\mu$ that appears in Brillinger's result. ii) Thus, $\|\hat\x-\mu\x_0\|_2^2$ is a natural measure of performance. iii) Estimation is possible even with $m<n$ measurements by taking advantage of the structure of $\x_0$. 

%


\subsection{Summary of Contributions}
Inspired by the work of Plan and Vershynin \cite{Ver}, and, motivated by recent advances on the precise analysis of the Generalized LASSO with linear measurements, this paper extends these latter results to the case of non-linear mesaurements. When the measurement matrix $\A$ has entries i.i.d. Gaussian (henceforth, we assume this to be the case without further reference), and the estimation performance is measured in a mean-squared-error sense, we are able to \emph{precisely} predict the asymptotic 
behavior of the error. The derived expression accurately captures the role of the link function $g$, the particular structure of $\x_0$, the role of the regularizer $f$, and, the value of the regularizer parameter $\la$. Further, it holds for all values of $\la$, and for a wide class of functions $f$ and $g$.

Interestingly, our result shows in a very precise manner that in large dimensions, modulo the information about the magnitude of $\x_0$, 
  the LASSO treats non-linear measurements exactly as if they were scaled and noisy linear measurements with scaling factor $\mu$ and  noise variance $\sigg$ defined as
 \begin{align}\label{eq:musig2}
\mu := \E[\gamma g(\gamma)], \quad\text{ and }\quad \sigg:=\E[(g(\gamma)-\mu\gamma)^2], \qquad\text{ for } \gamma\sim\Nn(0,1),
 \end{align}
where the expecation is with respect to both $\gamma$ and $g$.
In particular, when $g$ is such that $\mu\neq 0$\footnote{This excludes for example link functions $g$ that are even, but also see \cite[Sec.~2.2]{garnham2013note}}, then,
\begin{center}
\emph{the estimation performance of the Generalized LASSO with measurements of the form $y_i=g_i(\ab_i^T\x_0)$  is asymptotically the same as if the measurements were rather of the form $y_i=\mu \ab_i^T\x_0 + \sigma z_i$, with $\mu,\sigg$ as in \eqref{eq:musig2} and $z_i$ standard gaussian noise.}
\end{center}


\begin{figure*}[t!]
        \centering
        \includegraphics[width=0.6\textwidth]{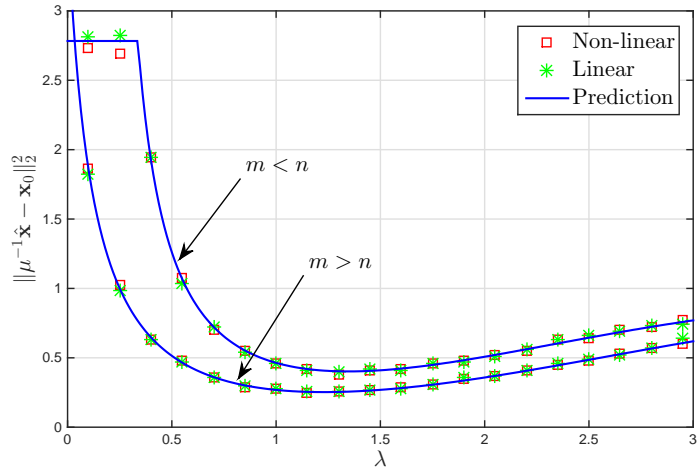}
        \caption{\footnotesize{Squared error of the $\ell_1$-regularized LASSO with non-linear measurements ($\square$) and with corresponding linear ones ($\star$) as a function of the regularizer parameter $\la$; both compared to the asymptotic prediction. Here, $g_i(x)=\text{sign}(x+0.3 z_i)$ with $z_i\sim\Nn(0,1)$. The unknown signal $\x_0$ is of dimension n=768 and has $\lceil{0.15n}\rceil$ non-zero entries (see Sec.~\ref{sec:sparse} for details). 
        The  different curves correspond to $\lceil0.75 n\rceil$ and $\lceil1.2 n\rceil$ number of measurements, respectively. Simulation points are averages over 20 problem realizations. }}
\label{fig:intro}
\end{figure*}

Recent analysis of the squared-error of the LASSO, when used to recover structured signals from noisy \emph{linear} observations, provides us with either precise predictions (e.g. \cite{ICASSP15,montanariLASSO}), or in other cases, with tight upper bounds (e.g. \cite{OTH13,DMM}). Owing to the established relation between non-linear and (corresponding) linear measurements, such results also characterize the performance of the LASSO in the presence of nonlinearities. We remark that some of the error formulae derived here in the  general context of non-linear measurements, have not been previously known even under the prism of linear measurements.



Figure \ref{fig:intro} serves as an illustration; the error with non-linear measurements matches well with the error of the corresponding linear ones and both are accurately predicted by our analytic expression.


Under the generic model in \eqref{eq:model}, which allows for $g$ to even be unspecified, $\x_0$ can, in principle, be estimated only up to a constant of proportionality \cite{Bri,li1989regression,Ver}. For example, if $g$ is uknown then any information about the norm $\|\x_0\|_2$ could be absorbed in the definition of $g$. The same is true when $g(x)=\text{sign}(x)$, eventhough $g$ might be known here. In these cases, what becomes important is the \emph{direction} of $\x_0$. Motivated by this, and, in order to simplify the presentation, we have assumed  throughout that \emph{$\x_0$ has unit Euclidean norm}\footnote{In \cite[Remark~1.8]{Ver}, they note that their results can be easily generalized to the case when $\|\x_0\|_2\neq 1$ by simply redifining $\bar g(x) = g(\|\x_0\|_2 x)$ and accordingly adjusting the values of the parameters $\mu$ and $\sigg$ in \eqref{eq:musig2}. The very same argument is also true in our case.     }, i.e. $\|\x_0\|_2=1$.
\subsection{Discussion of Relevant Literature}\label{sec:discuss}
\textbf{Extending an Old Result}. ~
Brillinger \cite{Bri} identified the asymptotic behavior of the estimation error of the LS solution $\hat\x_{LS}=(\A^T\A)^{-1}\A^T\y$ by showing that, when $n$ (the dimension of $\x_0$) is fixed,
 \begin{align}\label{eq:Bri}
 \lim_{{m}\rightarrow\infty} \sqrt{m} \|\hat\x_{LS} - \mu\x_0 \|_2 = {\sigma},
 \end{align}
 where  $\mu$ and $\sigg$ are same as in \eqref{eq:musig2}.  
 Our result can be viewed as a generalization  of the above in several directions. First, we extend \eqref{eq:Bri} to the regime where  $m/n=\delta\in(1,\infty)$ and both grow large by showing that
 \begin{align}\label{eq:Bri2}
 \lim_{n \rightarrow\infty} \|\hat\x_{LS} - \mu\x_0 \|_2 = \frac{\sigma}{\sqrt{\delta-1}}.
  \end{align} 
 Second, and most importantly, we consider solving the Generalized LASSO instead, to which LS is only a very special case. This allows versions of \eqref{eq:Bri2} where the error is finite even when  $\delta<1$ (e.g., see \eqref{eq:weVSVer}). 
 Note the additional challenges faced when considering the LASSO: i) $\hat\x$ no longer has a closed-form expression, ii) the result needs to additionally capture the role of $\x_0$, $f$, and, $\la$.

 \vspace{5pt}
\textbf{Motivated by Recent Work}.~
Plan and Vershynin consider a \emph{constrained} Generalized LASSO:
\begin{align}
\hat\x_{\text{C-LASSO}} = \arg\min_{\x\in\Kc}{\|\y-\A\x\|_2},
\end{align}
with $\y$ as in \eqref{eq:model} and $\Kc\subset\R^n$ some known set (not necessarily convex).
 In its simplest form, their result shows that when $m\gtrsim \DK$ then with high probability,
\begin{align}\label{eq:Ver}
\|\hat\x_{\text{C-LASSO}} - \mu\x_0\|_2\lesssim \frac{\sigma\sqrt{\DK}+\zeta}{\sqrt{m}}.
\end{align} 
Here, $\DK$ is the Gaussian width, a specific measure of complexity of the constrained set $\Kc$ when viewed from $\mu\x_0$. For our purposes, 
it suffices to remark that if $\Kc$ is properly chosen, and, if $\mu\x_0$ is on the boundary of $\Kc$, then $\DK$ is less than  $n$. Thus, estimation is in principle is possible with $m<n$ measurements. The parameters $\mu$ and $\sigma$ that appear in \eqref{eq:Ver} are the same as in \eqref{eq:musig2} and $\zeta:=\E[(g(\gamma)-\mu\gamma)^2\gamma^2]$.  Observe that, in contrast to \eqref{eq:Bri} and to the setting of this paper, the result  in \eqref{eq:Ver} is non-asymptotic. Also, it suggests the critical role played by $\mu$ and $\sigma$. On the other hand, \eqref{eq:Ver} is only an upper bound on the error, and also, it suffers from unknown absolute proportionality constants (hidden in $\lesssim$).

Moving the analysis into an asymptotic setting, our work expands upon the result of \cite{Ver}. First, we consider the regularized LASSO instead, which is more commonly used in practice. Most importantly, we improve the loose upper bounds into \emph{precise} expressions. In turn, this proves in an exact manner the role played by $\mu$ and $\sigg$ to which \eqref{eq:Ver} is only indicative. For a direct comparison with \eqref{eq:Ver} we mention the following result which follows from our analysis (we omit the proof for brevity). Assume $\Kc$ is convex, $m/n=\delta\in(0,\infty)$, $\DK/n=\rho\in (0,1]$ and $n\rightarrow\infty$. Also, $\delta>\rho$. Then, \eqref{eq:Ver} yields an upper bound $C\sigma\sqrt{\rho/\delta}$ to the error,  for some constant $C>0$. Instead, we  show  
\begin{align}\label{eq:weVSVer}
\|\hat\x_{\text{C-LASSO}} - \mu\x_0\|_2 \leq \sigma\frac{\sqrt{\rho}}{\sqrt{\delta-\rho}}.
\end{align}
\vspace{5pt}
\textbf{Precise Analysis of the LASSO With Linear Measurements}.~
The first \emph{precise} formula predicting the limiting behavior of the LASSO reconstruction error wer established in \cite{DMM,montanariLASSO,StoLASSO}. The authors of \cite{DMM,montanariLASSO} consider  the $\ell_2^2$-LASSO with $\ell_1$-regularization and the analysis is based on the  the Approximate Message Passing
(AMP) framework\cite{AMP}; also  \cite{Armeen,malekiComplexLASSO} for extensions. A more general line of work, \cite{StoLASSO,OTH13,ICASSP15,ISIT15_ell22} studies the problem using a recently developed framework that is based on Gordon's Gaussian min-max Theorem (GMT) \cite{StoLASSO,COLT15}. The GMT framework was initially used by Stojnic \cite{StoLASSO} to derive tight upper bounds on the constrained LASSO with $\ell_1$-regularization; \cite{OTH13} generalized those to general convex regularizers and also to the $\ell_2$-LASSO; the case of the $\ell_2^2$-LASSO was studied in \cite{ISIT15_ell22}. Those bounds hold for all values of SNR, but they become tight only in the high-SNR regime. A precise error expression was derived in \cite{ICASSP15} for the $\ell_2$-LASSO with $\ell_1$-regularization under a gaussianity assumption on the distribution of the non-zero entries of $\x_0$. When measurements are linear, our Theorem \ref{coro:sparse} generalizes this assumption; in its current form, it provides the first-known counterpart of the main result of \cite{montanariLASSO} for the $\ell_2$-LASSO. Our main Theorem \ref{thm:main} provides error predictions for regularizers going beyond the $\ell_1$-norm, e.g. $\ell_{1,2}$-norm, nuclear norm, which appear to be novel.
When it comes to non-linear measurements, to the best of our knowledge, this paper is the first to derive asymptotically precise results on the performance of any LASSO-type algorithm.


%






\section{Results}
\subsection{Modeling Assumptions} \label{sec:model}
\textbf{{Unknown structured signal}}. We let $\x_0\in\R^n$ represent the unknown signal vector. We  assume that
$$\x_0={\xx_0}/{\|\xx_0\|_2},$$
with $\xx_0$ sampled from a probability density $p_{\overline{\x}_0}$ in $\R^n$. Thus, $\x_0$ is deterministically of unit Euclidean-norm (this is mostly to simplify the presentation, see Footnote 4). Information about the structure of  $\xx_0$ (and correspondingly of $\x_0$) is encoded in $p_{\xx_0}$. For instance, to study an $\xx_0$ which is sparse, it is typical to assume that its entries are i.i.d. $\xx_{0,i}\sim(1-\rho)\delta_0+\rho q_{\XX_0}$, where $\rho\in(0,1)$ becomes the normalized sparsity level, $q_{\XX_0}$ is a scalar p.d.f. and $\delta_0$ is the Dirac delta function\footnote{Such models in place for studying structured signals have been widely used in the relevant literature, e.g. \cite{DJ94,DMM,donoho2013accurate}. In fact, the results here continue to hold as long as the marginal distribution of $\xx_0$ converges to a given distribution (as in \cite{malekiComplexLASSO,montanariLASSO}).}. 


%
%
%

\vspace{5pt}
\textbf{{Regularizer}}. We consider \emph{convex} regularizers $f:\R^n\rightarrow\R$. 

\vspace{5pt}
\textbf{Measurement matrix}. The entries of $\A\in\R^{m\times n}$ are i.i.d. $\Nn(0,1)$.

\vspace{5pt}
\textbf{Measurements and Link-function}.  We observe $\y=\gv(\A\x_0)$ where $\gv$ is a (possibly random) map from $\R^m$ to $\R^m$ and $\gv(\ub)=[g_1(u_1),\ldots,g_m(u_m)]^T$. Each $g_i$ is  i.i.d. from a real valued random function $g$ for which $\mu$ and $\sigg$ are defined in \eqref{eq:musig2}.  We assume that $\mu$ and $\sigg$ are nonzero and bounded.


\vspace{5pt}
\textbf{Asymptotics.} 
 We study a linear asymptotic regime. In particular, we consider a 
sequence of problem instances $\{\xx_0^\mn, \A^\mn, f^\mn, m^\mn\}_{n\in\mathbb{N}}$ indexed by $n$ such that 
 $\A^\mn\in\R^{m\times n}$ has entries i.i.d. $\Nn(0,1)$, $f^\mn:\R^n\rightarrow\R$ is proper convex, and,
$m:=m^\mn$ with $m= {\delta n}, \delta\in(0,\infty)$. 
We further require that the following conditions hold:
\begin{enumerate}[(a)]

\item $\xx_0^\mn$ is sampled from a probability density $p^\mn_{\overline{\x}_0}$ in $\R^n$ with one-dimensional marginals that are independent of $n$ and have bounded second moments. Furthermore, 
$n^{-1}\|\overline{\x}^\mn_0\|_2^2\rP \sigma_x^2=1$.


\item For any $n\in\mathbb{N}$ and any $\|\x\|_2\leq C$, it holds $n^{-1/2}f(\x)\leq c_1$ and $n^{-1/2}\max_{\s\in\partial{f}^\mn(\x)}\|\s\|_2\leq c_2$, for constants $c_1,c_2,C\geq0$ independent of $n$. 

\end{enumerate}

{In (a), we used ``$\rP$" to denote convergence in probability as $n\rightarrow\infty$. The assumption $\sigma_x^2=1$ holds without loss of generality, and, is only necessary to simplify the presentation. In (b), $\partial f(\x)$ denotes the subdifferential of $f$ at $\x$. The condition itself is no more than a normalization condition on $f$.}


Every such sequence $\{\xx_0^\mn, \A^\mn, f^\mn \}_{n\in\mathbb{N}}$ generates a sequence $\{\x_0^\mn, \y^\mn\}_{n\in\mathbb{N}}$ where $\x_0^\mn := {\xx^\mn_0}/{\|\xx^\mn_0\|_2}$ and $\y^\mn := \gv^\mn(\A\x_0)$.  @hen clear from the context, we  drop the superscript $(n)$.



\subsection{General Result}
Let $\{\x_0^\mn, \A^\mn, f^\mn, \y^\mn \}_{n\in\mathbb{N}}$ be a sequence of problem instances that satisfies the conditions of Section \ref{sec:model}. With these, define the sequence 
$\{\hat\x^\mn \}_{n\in\mathbb{N}}$
of solutions to the corresponding LASSO problems for fixed $\la>0$:
\begin{align}\label{eq:genLASSO}
\hat\x^\mn := \min_{\x} \frac{1}{\sqrt{n}} \left\{ \|\y^\mn - \A^\mn\x\|_2 + \la f^\mn(\x) \right\}.
\end{align}

The main contribution of this paper is a precise evaluation of  $\lim_{n\rightarrow\infty}\|\mu^{-1}\x^\mn-\x^\mn_0\|_2^2$ with high probability over the randomness of $\A$, of $\x_0$, and of  $g$.
To state the result in a general framework, we require a further assumption on $p_{\xx_0}^\mn$ and $f^\mn$. Later in this section we illustrate how this assumption can be naturally met.  We write $f^*$ for the Fenchel's conjugate of $f$, i.e., $f^*(\vb):=\sup_{\x}\x^T\vb-f(\x)$; also, we call the proximal function of $f$ to be  $\prox_{f,\tau}(\vb) := \min_{\x}\{\frac{1}{2}\|\vb-\x\|_2^2 + \tau f(\x)\}$.

\vspace{5pt}
\begin{ass}\label{ass:tech}
We say Assumption \ref{ass:tech} holds 
if  for all non-negative  constants $c_1,c_2,c_3\in\R$  the point-wise limit of $\frac{1}{n}\prox_{\sqrt{n}(f^*)^\mn,c_3}\left(c_1\h + c_2 \xx_0\right)$ exists with probability one over $\h\sim\Nn(0,\I_n)$ and $\xx_0\sim p^\mn_{\xx_0}$. Then, we denote the limiting value as $F(c_1,c_2,c_3)$.
%



\end{ass}
\vspace{8pt}
\begin{thm}[Non-linear$=$Linear]\label{coro:main} 
Consider the asymptotic setup of Section \ref{sec:model} and let Assumption \ref{ass:tech} hold. Recall $\mu$ and $\sigma^2$ as in \eqref{eq:musig2} and let $\hat\x$ be the minimizer of the Generalized LASSO in \eqref{eq:genLASSO} for fixed $\la>0$ and for measurements given by \eqref{eq:model}.
Further let  $\hat\x^\text{\lin}$ be the solution to the Generalized LASSO when used with \emph{linear} measurements of the form
$\y^\text{\lin} = \A(\mu\x_0) + \sigma\z$,
where $\z$ has entries i.i.d. standard normal.
Then, in the limit of $n\rightarrow\infty$, with probability one,
\begin{align*}
\|\hat\x-\mu\x_0\|_2^2 = \|\hat\x^\text{lin}-\mu\x_0\|_2^2.
\end{align*}
\end{thm}

\vspace{3pt}
Theorem \ref{coro:main} relates in a very precise manner the error of the Generalized LASSO under  non-linear measurements to the error of the same algorithm when used  under appropriately scaled noisy linear measurements. Theorem \ref{thm:main} below, derives an asymptotically exact expression for the error.
\vspace{8pt}
\begin{thm}[Precise Error Formula]\label{thm:main}
Under the same assumptions of Theorem \ref{coro:main} and $\delta:=m/n$, it holds,  with probability one,
$$
\lim_{n\rightarrow\infty}\|\hat\x-\mu\x_0\|_2^2 = \alpha_*^2,
$$
where $\alpha_*$ is the unique optimal solution to the convex program
\begin{align}\label{eq:maxmin}
\max_{\substack{0\leq\beta\leq 1 \\ \tau\geq 0}}\min_{\alpha\geq 0} \beta\sqrt{\delta}\sqrt{\alpha^2+\sigma^2} - \frac{\alpha\tau}{2} + \frac{\mu^2\tau}{2\alpha}-\frac{\alpha\la^2}{\tau}F\left(\frac{\beta}{\la},\frac{\mu\tau}{\la\alpha},\frac{\tau}{\la\alpha}\right).
\end{align}
Also, the optimal cost of the LASSO  in \eqref{eq:genLASSO} converges to the optimal cost of the program in \eqref{eq:maxmin}.
 \end{thm}

\vspace{5pt}
Under the stated conditions, Theorem \ref{thm:main} proves that the limit of $\|\hat\x-\mu\x_0\|_2$ exists and is equal to the \emph{unique} solution of the optimization program in \eqref{eq:maxmin}. Notice that this is a \emph{deterministic} and \emph{convex} optimization, which only involves three \emph{scalar} optimization variables. Thus, the optimal $\alpha_*$ can, in principle, be efficiently numerically computed. In many specific cases of interest, with some extra effort,  it is possible to yield simpler expressions for $\alpha_*$, 
e.g. see Theorem \ref{coro:sparse} below. The role of the normalized number of measurement $\delta=m/n$, of the regularizer parameter $\la$, and, that of $g$, through $\mu$ and $\sigg$, are explicit in \eqref{eq:maxmin}; the structure of $\x_0$ and the choice of the regularizer $f$ are implicit in $F$.
Figures \ref{fig:intro}-\ref{fig:main} illustrate the accuracy of the prediction of the theorem in a number of different settings. 
The proofs of both the Theorems are deferred  to Appendix  \ref{sec:main_app}.
In the next sections, we specialize Theorem \ref{thm:main} to the cases of sparse, group-sparse and low-rank signal recovery. 


\subsection{Examples}
\subsubsection{Sparse Recovery}\label{sec:sparse}

Assume each entry $\xx_{0,i}, i=1,\ldots,n$ is sampled i.i.d. from a distribution 
\begin{align}
p_{\XX_0}(x) = (1-\rho)\cdot\delta_0(x) + \rho\cdot q_{\XX_0}(x),\label{eq:sparseDist}
\end{align}
 where $\delta_0$ is the delta Dirac function, $\rho\in(0,1)$ and $q_{\XX_0}$ a probability density function with second moment normalized to $1/\rho$ so that  condition (a) of Section \ref{sec:model} is satisfied. Then, $\x_0=\xx_0/\|\xx_0\|_2$ is $\rho n$-sparse on average and has unit Euclidean norm. 
 Letting $f(\x)=\|\x\|_1$ also satisfies condition (b). Let us now check Assumption 1. 
The Fenchel's conjugate of the $\ell_1$-norm is simply the indicator function of the $\ell_\infty$ unit ball. Hence, 
without much effort,
\begin{align}
\frac{1}{n}\prox_{\sqrt{n}(f^*)^\mn,c_3}\left(c_1\h + c_2 \xx_0\right) &= \frac{1}{2n} \sum_{i=1}^{n}\min_{|v_i|\leq 1} ( v_i - (c_1 \h_i + c_2 \xx_{0,i} ))^2\nn \\&= \frac{1}{2n} \sum_{i=1}^{n} \eta^2(c_1 \h_i + c_2 \xx_{0,i}; 1),\label{eq:eta_sparse}
\end{align}
where we have denoted 
\begin{align}\label{eq:soft}
\eta(x;\tau):=({x}/{|x|})\left(|x|-\tau\right)_+
\end{align}
 for the soft thresholding operator. An application of the weak law of large numbers to see that the limit of the expression in \eqref{eq:eta_sparse} equals $F(c_1,c_2,c_3):=\frac{1}{2}\E\left[\eta^2(c_1 h + c_2 \XX_0; 1)\right],$
where the expectation is over $h\sim\Nn(0,1)$ and $\XX_0\sim p_{\XX_0}$. With all these, Theorem \ref{thm:main} is applicable.
We have put  extra effort in order to obtain the following equivalent but more insightful characterization of the error, as stated below and proved in the Appendix.

\begin{figure*}[t!]
    \centering
    \begin{subfigure}[t]{0.5\textwidth}
        \centering
        \includegraphics[width=\textwidth]{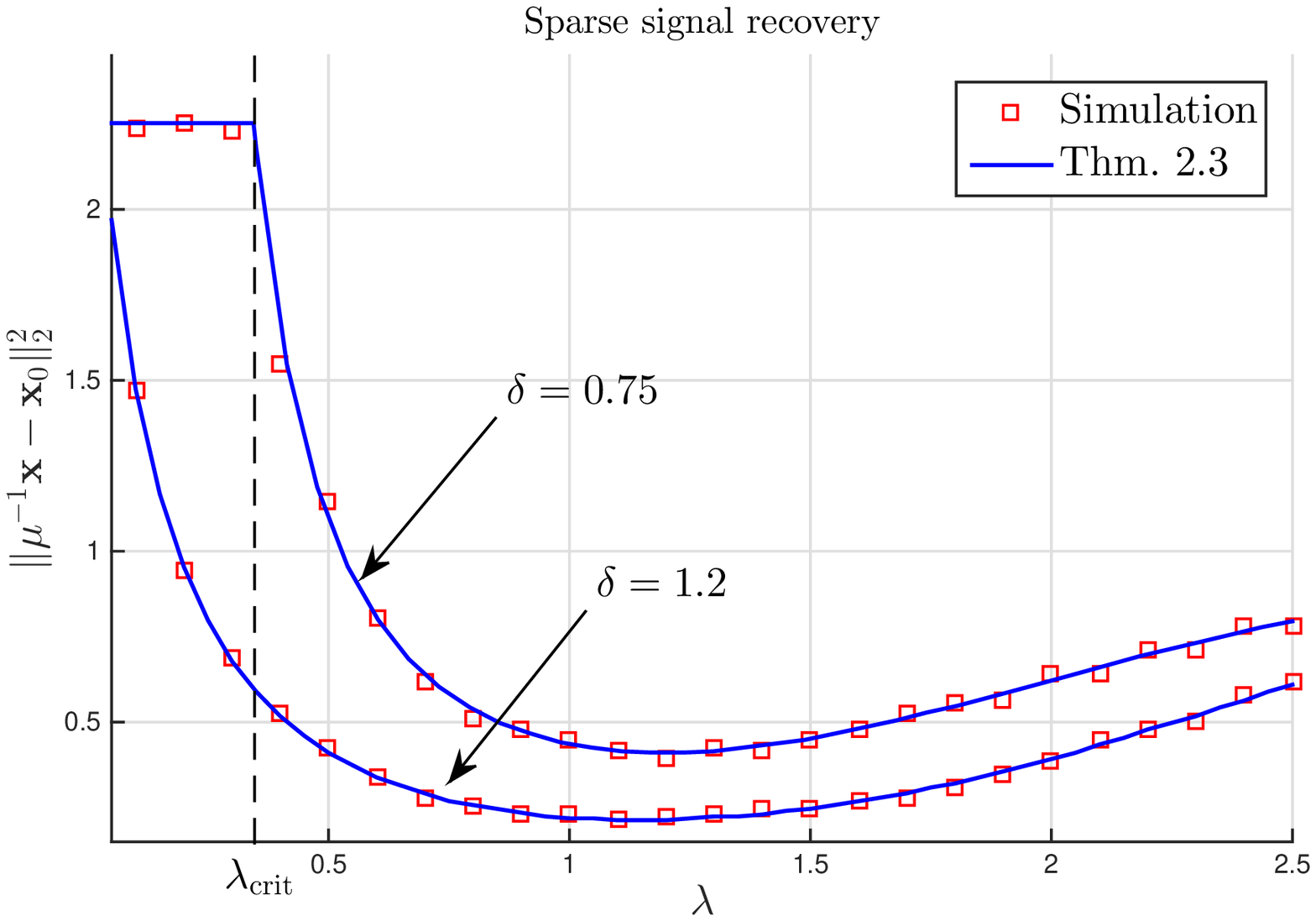}\label{fig:sparse}

    \end{subfigure}%
    ~ 
    \begin{subfigure}[t]{0.5\textwidth}
        \centering
        \includegraphics[width=\textwidth]{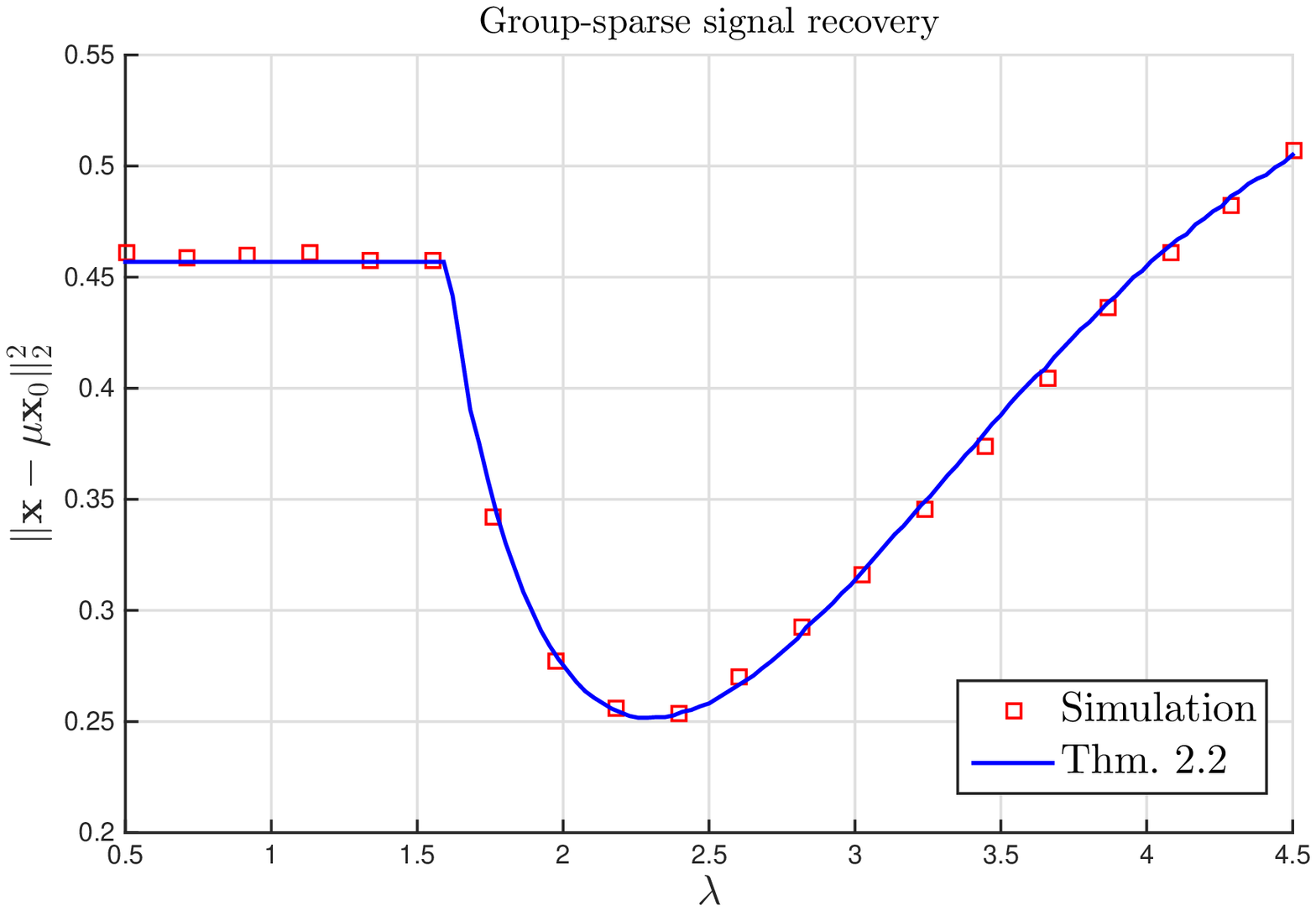}
                \label{fig:gs}
    \end{subfigure}
    \caption{\small{ Squared error of the LASSO as a function of the regularizer parameter compared to the asymptotic predictions. Simulation points represent averages over 20 realizations. (a) Illustration of Thm.~\ref{coro:sparse} for $g(x)=\text{sign}(x)$, $n=512$, $p_{\XX_0}(+1)=p_{\XX_0}(+1)=0.05$, $p_{\XX_0}(+1)=0.9$ and two values of $\delta$, namely $0.75$ and $1.2$. (b)  Illustration of Thm.~\ref{thm:main} for $\x_0$ being group-sparse as in Section \ref{sec:gs} and $g_i(x) = \mathrm{sign}(x+0.3z_i)$. In particular, $\x_0$ is composed of $t=512$ blocks of block size $b=3$. Each block is zero with probability $0.95$, otherwise its entries are i.i.d. $\Nn(0,1)$. Finally, $\delta=0.75$.}}\label{fig:main}
\end{figure*}

\begin{thm}[Sparse Recovery] \label{coro:sparse}
If $\delta>1$, then define $\lac=0$. Otherwise, let $\lac, \kac$  be the unique pair of solutions to the following set of equations:
\begin{numcases}{}
\kappa^2\delta = \sigma^2 + \E\left[( \eta(\kappa h + \mu \XX_0; \kappa\la) - \mu\XX_0 )^2\right],\label{eq:tau}\\
\kappa\delta = \E[( \eta(\kappa h + \mu \XX_0; \kappa\la)\cdot h )],\label{eq:beta}
\end{numcases}
where $h\sim\Nn(0,1)$ and is independent of $\XX_0\sim p_{\XX_0}$.
Then, for any $\la>0$, with probability one,
$$
\lim_{n\rightarrow\infty}\|\hat\x-\mu\x_0\|_2^2 = \begin{cases} \delta\kac^2 - \sigma^2 &,\la\leq\lac, \\ \delta\kappa_*^2(\la) - \sigma^2&, \la\geq\lac, \end{cases}
$$
where $\kappa^2_*(\la)$ is the unique solution to \eqref{eq:tau}.
%
\end{thm}


Figures \ref{fig:intro} and \ref{fig:main}(a) validate the prediction of the theorem, for different signal distributions, namely $q_{\XX_0}$ being Gaussian and Bernoulli, respectively. For the case of  compressed ($\delta<1$) measurements, observe the two different regimes of operation, one for $\la\leq\lac$ and the other for $\la\geq\lac$, precisely as they are predicted by the theorem (see also \cite[Sec.~8]{OTH13}). 
The special case of Theorem \ref{coro:sparse} for which $q_{\XX_0}$ is Gaussian has been previously studied in \cite{ICASSP15}. Otherwise, to the best of our knowledge, this is the first precise analysis result for the $\ell_2$-LASSO stated in that generality. Analogous result, but via different analysis tools, has only been known for the $\ell_2^2$-LASSO as appears in \cite{montanariLASSO}.

\vspace{8pt}
\subsubsection{Group-Sparse Recovery}\label{sec:gs}
Let $\xx_0\in\R^n$ be composed of $t$ non-overlapping blocks of constant size $b$ each such that $n=t\cdot b$. Each block $[\xx_0]_i,i=1,\ldots,t$ is sampled i.i.d. from a probability density in $\R^b$: $p_{\XX_0}(\x) = (1-\rho)\cdot\delta_0(\x) + \rho\cdot q_{\XX_0}(\x), \x\in\R^b$, where $\rho\in(0,1)$. Thus, $\xx_0$ is a $\rho t$-block-sparse on average. We operate in the regime of linear measurements $m/n=\delta\in(0,\infty)$. As is common we use the $\ell_{1,2}$-norm to induce block-sparsity, i.e., $f(\x)=\sum_{i=1}^{t}\| [\xx_0]_i \|_2$; with this, \eqref{eq:genLASSO} is often referred to as group-LASSO in the literature \cite{groupLASSO}. It is not hard to show that Assumption 1 holds with $F(c_1,c_2,c_3):=\frac{1}{2b}\E\left[\|\etav(c_1 \h + c_2 \XX_0; 1)\|_2^2\right],$ where  $\etav(\x;\tau) =  {\x}/{\|\x\|}\left(\|\x\|_2-\tau\right)_+, \x\in\R^b$ is the vector soft thresholding operator and $h\sim\Nn(0,\I_b)$, $\XX_0\sim p_{\XX_0}$ and are independent. Thus Theorem \ref{thm:main} is applicable in this setting; Figure \ref{fig:main}(b) illustrates the accuracy of the prediction.


\vspace{8pt}
\subsubsection{Low-rank Matrix Recovery}

Let ${\XXb_0}\in\R^{d\times d}$ be an unknown matrix of rank $r$, in which case, $\xx_0=\mathrm{vec}(\XXb_0)$ with  $n=d^2$. Assume $m/d^2=\delta\in(0,\infty)$ and $r/d=\rho\in(0,1)$. As usual in this setting, we consider nuclear-norm regularization; in particular, we choose $f(\x)=\sqrt{d}\|\X\|_*$. Each subgradient $\Sb\in\partial f(\X)$ then satisfies $\|\Sb\|_F \leq d$ in agreement with assumption (b) of Section \ref{sec:model}. Furthermore, for this choice of regularizer, we have
\begin{align}
&\frac{1}{n}\prox_{\sqrt{n}(f^*)^\mn,c_3}\left(c_1\Hb + c_2 \XXb_0\right) = \frac{1}{2d^2} \min_{\|\Vb\|_2\leq \sqrt{d}} \| \Vb - (c_1 \Hb + c_2 \XXb_{0} )\|_F^2\nn \\&= \frac{1}{2d} \min_{\|\Vb\|_2\leq 1} \| \Vb - d^{-1/2}(c_1 \Hb + c_2 \XXb_{0} )\|_F^2
 = \frac{1}{2d} \sum_{i=1}^{d} \eta^2\left( \mathrm{s_i}\left( d^{-1/2}(c_1 \Hb + c_2 \XXb_{0} ) \right) ; 1\right) ,\nn
\end{align}
where $\eta(\cdot;\cdot)$ is as in \eqref{eq:soft}, $\mathrm{s_i}(\cdot)$ denotes the $i^\text{th}$ singular value of its argument and $\Hb\in\R^{d\times d}$ has entries $\Nn(0,1)$. If conditions are met such that the empirical distribution of the singular values of (the sequence of random matrices) $c_1\Hb+c_2\XXb_0$ converges asymptotically to a limiting distribution, say $q(c_1,c_2)$, then $F(c_1,c_2,c_3):=\frac{1}{2}\E_{x\sim q(c_1,c_2)}\left[\eta^2( x ; 1)\right],$ and Theorem \ref{coro:main}--\ref{thm:main} apply. For instance, this will be the case if $d^{-1/2}\XXb_0=\Ub\Sb\Vb^t$, where $\Ub,\Vb$ unitary matrices and $\Sb$ is a diagonal matrix whose entries have a given marginal distribution with bounded moments (in particular, independent of $d$). We leave the details and the problem of (numerically) evaluating $F$ for future work.

\subsection{An Application to $q$-bit Compressive Sensing}\label{sec:qbit}
\subsubsection{Setup}

Consider recovering a \emph{sparse}  unknown signal $\x_0\in\R^n$ from scalar q-bit quantized linear measurements.
 Let  $\tb:=\{t_0=0, t_1,\ldots, t_{L-1}, t_L=+\infty\}$ represent a (symmetric with respect to $0$) set of decision thresholds and $\ellb:=\{\pm \ell_1,\pm \ell_2,\ldots,\pm \ell_L\}$ the corresponding representation points, such that $L=2^{q-1}$. Then, quantization of a real number $x$ into $q$-bits can be represented 
as
\begin{align}
\Qc_q(x,\ellb,\tb) = \sign(x) \sum_{i=1}^{L}\ell_i\one_{\{ t_{i-1}\leq |x| \leq t_{i} \}},\nn
\end{align}
where $\one_{\Sc}$ is the indicator function of a set $\Sc$. For example, 1-bit quantization with level $\ell$ corresponds to $Q_1(x,\ell) = \ell\cdot\sign(x)$. 
The measurement vector $\y=[y_1,y_2\ldots,y_m]^T$ takes the form
\begin{align}\label{eq:Qmeas}
y_i = \Qc_q(\ab_i^T\x_0,\ellb,\tb),\quad i=1,2,\ldots,m,
\end{align}
where $\ab_i^T$'s are the rows of a measurement matrix $\A\in\R^{m\times n}$, which is henceforth assumed i.i.d. standard Gaussian. 
We use the LASSO to obtain an estimate $\hat\x$ of $\x_0$ as
\begin{align}\label{eq:QLASSO}
\hat\x:=\arg\min_\x\|\y-\A\x\|_2 + \la \|\x\|_1.
\end{align} 
Henceforth, we assume for simplicity that $\|\x_0\|_2=1$. Also, in our case, $\mu$ is known since $g=Q_q$ is known; thus, is reasonable to scale the solution of \eqref{eq:QLASSO} as $\mu^{-1}\hat\x$  and consider the error quantity $\|\mu^{-1}\hat\x-\x_0\|_2$ as a measure of estimation performance.
%
Clearly, the error depends (besides others)  on the number of bits $q$, on the choice of the decision thresholds $\tb$ and on the quantization levels $\ellb$. An interesting question of practical importance becomes how to optimally choose these to achieve less error. As a running example for this section, we seek optimal quantization thresholds and corresponding levels
\begin{align}\label{eq:Qobj}
(\tb_*,\ellb_*) = \arg\min_{\tb,\ellb} \|\mu^{-1}\hat\x-\x_0\|_2,
\end{align}
while keeping all other parameters such as the number of bits $q$ and  of measurements $m$  fixed.

\subsubsection{Consequences of Precise Error Prediction}\label{sec:conc}
Theorem \ref{coro:main} shows  that $
\|\mu^{-1}\hat\x-\x_0\|_2 = \|\hat\x^\lin-\x_0\|_2,
$ where $\hat\x^{\lin}$ is the solution to \eqref{eq:QLASSO}, but only, this time with a measurement vector $\y^{\lin}=\A\x_0 + \frac{\sigma}{\mu}\z$, where $\mu,\sigma$ as in \eqref{eq:Qmusig2}  and $\z$ has entries i.i.d. standard normal. Thus, lower values of the ration $\sigg/\mu^2$ correspond to lower values of the error and the design problem posed in \eqref{eq:Qobj} is equivalent to the following simplified one:
\begin{align}\label{eq:Qobj2}
(\tb_*,\ellb_*) = \arg\min_{\tb,\ellb} \frac{\sigg(\tb,\ellb)}{\mu^2(\tb,\ellb)}.
\end{align}
To be explicit, $\mu$ and $\sigg$ above can be easily expressed from \eqref{eq:musig2} after setting $g=\Qc_q$ as follows:
\begin{align}\label{eq:Qmusig2}
&\mu:=\mu(\ellb,\tb) = \sqrt{\frac{2}{\pi}}\sum_{i=1}^{L}\ell_i\cdot\left(e^{-t_{i-1}^2/2}-e^{-t_{i}^2/2}\right)
\quad\text{and}\quad\sigg:=\sigg(\ellb,\tb):=\tauu-\mu^2,\\ &\text{where,}\quad \tauu:=\tauu(\ellb,\tb) = 2\sum_{i=1}^{L}\ell_i^2\cdot\left(Q(t_{i-1})-Q(t_{i}\right)) \quad\text{and}\quad Q(x) = \frac{1}{\sqrt{2\pi}}\int_{x}^\infty\exp(-u^2/2)\mathrm{d}u.\nn
\end{align}

\subsubsection{An Algorithm for Finding Optimal Quantization Levels and Thresholds}\label{sec:claim}
In contrast to the initial problem in \eqref{eq:Qobj}, the optimization involved in \eqref{eq:Qobj2} is explicit in terms of the variables $\ellb$ and $\tb$, but, is still hard to solve in general. Interestingly, we show in the Appendix .
that the popular Lloyd-Max (LM) algorithm can be an effective algorithm for solving \eqref{eq:Qobj2}, since the values to which it converges are stationary points of the objective in \eqref{eq:Qobj2}. Note that this is not a directly obvious result since the classical objective of the LM algorithm is minimizing the quantity $\E[\|\y - \A\x_0\|_2^2]$ rather than $\E[\|\mu^{-1}\hat\x-\x_0\|_2^2]$. 


%
%

%

\bibliography{./../compbib}

\newpage
\appendix[Proof of Theorems \ref{coro:main} \& \ref{thm:main}]\label{sec:main_app}
\subsection{Theorem \ref{thm:main}}
We start with the proof of Theorem \ref{thm:main}. Theorem \ref{coro:main} will follow as a direct corollary of this result.

Assume a sequence of problem instances as described in Section \ref{sec:model}. To keep notation simple,  we simply use $\|\vb\|$ (rather than $\|\vb\|_2$) for the Euclidean norm of $\vb$ and we shall also drop the superscript $\mn$ when referring to elements of the sequence.  Thus, we write
\beq\label{Eq1}
\hat{\x}=\arg\min_{\x} \frac{1}{\sqrt{n}}\|\vec{g}(\A\x_0)-\A\x\|+\frac{\la}{\sqrt{n}}\f(\x),
\eeq
but it is to be understood that the above actually produces a sequence of solutions $\hat\x^\mn$ indexed by $n$. Our goal is to characterize the nontrivial limiting behavior of $\|\hat\x-\mu\x_0\|_2$. 

We start with a simple but useful change of variables $\w:=\x-\mu\x_0$, to directly have a handle on the error vector $\w$. Then, \eqref{Eq1} becomes:
\begin{align}
\hat{\w}&:=\arg\min_{\w} \frac{1}{\sqrt{n}}\|\vec{g}(\A\x_0)- \mu\A\x_0 - \A\w\|+\frac{\la}{\sqrt{n}}\f(\mu\x_0+\w) \nn \\ 
&=
\arg\min_{\w}\max_{\|\ub\|\leq 1} \frac{1}{\sqrt{n}}	(-\ub^T\A\w +\ub^T(\vec{g}(\A\x_0)-\mu \A\x_0))	+\frac{\la}{\sqrt{n}}	\f(\mu\x_0+\w),\label{eq2}
\end{align}
where the second line follows after using the fact $\|\vb\|=\max_{\|\ub\|_2\leq 1}\ub^T\vb$.

\subsubsection{A Key Decomposition}
The first key step in the proof is a trick adapted from the proofs of \cite[Lem.~4.3]{Ver} and \cite[Thm.~1.3]{plan2014high}. Until further notice, we condition on $\x_0$. Also, we repeatedly make use of the assumption that $\|\x_0\|=1$ without direct reference. 
The trick amounts to decomposing each measurement vector $\ab_i$ in its projection on the direction of $\x_0$ and its orthogonal complement. Denoting $\Pp=(\I-\x_0\x_0^T)$ for the projector onto the orthogonal complement of the span of $\x_0$ (recall $\|\x_0\|_2=1$), we have 
$
\ab_i^T = (\ab_i^T\x_0)\x_0^T + \ab_i^T\Pp,
$ or, in matrix form:
$$
\A = (\A\x_0)\x_0^T + \A\Pp.
$$
Then, \eqref{eq2} becomes:
\begin{align}
\min_{\w}\max_{\|\ub\|\leq 1} \frac{1}{\sqrt{n}}	-\ub^T\A\Pp\w +\ub^T(\vec{g}(\A\x_0)-\mu \A\x_0 - (\A\x_0)\x_0^T\w)	+\frac{\la}{\sqrt{n}}	\f(\mu\x_0+\w).\label{eq3}
\end{align}
Now, we use the Gaussianity assumption on the entries of $\A$ to see that $\A\Pp$ is \emph{independent} of $\A\x_0$. It can then be shown (see \cite[pg.~13]{Ver}) that $\A\Pp$ is also \emph{independent} of $(\vec{g}(\A\x_0)-\mu \A\x_0)$; thus, $\A\Pp\w$ is independent of the rest terms in in \eqref{eq3}. This shows that the objective function of \eqref{eq3} is distributed identically even after replacing the $\A\Pp\w$ with $\G\Pp\w$, where $\G$ is an independent copy of $\A$. After all these, \eqref{eq3} is identically distributed with the following:
\begin{align}
\min_{\w}\max_{\|\ub\|\leq 1} \frac{1}{\sqrt{n}}\{ -\ub^T\G\Pp\w +\ub^T(\z_\e - (\x_0^T\w)\e)	\}+\frac{\la}{\sqrt{n}}	\f(\mu\x_0+\w).\label{eq4}
\end{align}
where $\G$ and  $\e:=\A\x_0$ have entries i.i.d. standard normal and are independent of each other. Also, $\z_\e:=\vec{g}(\e)-\mu\e$ for convenience.

\subsubsection{Applying the cGMT}
After the decomposition step in the previous section, we have transformed the initial problem to that of analyzing the (probabilistically) equivalent one in \eqref{eq4}. In particular, we wish to evaluate the limiting behavior of $\|\hat\w\|_2$, i.e. the norm of the minimizer of the optimization in \eqref{eq4}. The analysis is possible thanks to the convex Gaussian Min-max Theorem (cGMT) \cite[Thm.~1]{COLT15}, which is a stronger version of the classical result of Gordon \cite{gorThm} in the presence of additional convexity assumptions. According to the cGMT, the analysis of a Primary Optimization (PO) problem that is of the form
\begin{align}\label{eq:PO}
\min_{\vb\in\Sc_{\vb}}\max_{\ub\in\Sc_{\ub}} \ub^T\G\vb + \psi(\vb,\ub),
\end{align}
with $\G$ being i.i.d. Gaussian, $\Sc_\vb,\Sc_\ub$ convex, compact sets and $\psi$ a convex-concave function, can be carried out via analyzing a corresponding Auxiliary Optimization problem (AO), which is defined as
\begin{align}\label{eq:AO}
\min_{\vb\in\Sc_{\vb}}\max_{\ub\in\Sc_{\ub}} \|\vb\|\g^T\ub + \|\ub\|\h^T\vb + \psi(\vb,\ub).
\end{align}
In \eqref{eq:AO}, $\g$ and $\h$ are i.i.d. standard Gaussian vectors of appropriate size. To apply the theorem, identify $\vb:=\Pp\w$ in \eqref{eq4} and the appearance of the bilinear term $\ub^T\G\vb$ as in \eqref{eq:PO}. Also, the rest of the objective function in \eqref{eq4} is convex in $\Pp\w$ (where we have used the convexity of $f$) and linear (thus, concave) in $\ub$. Overall, \eqref{eq4} is in the appropriate format of a (PO) problem as in \eqref{eq:PO}. The only technical caveat is that the minimization over $\w$ in it appears unconstrained. For this, we assume that the minimizer of \eqref{eq4} satisfies $\|\hat\w\|\leq K_\w$ for sufficiently large constant $K_\w>0$ \emph{independent of $n$}. If our assumption is valid, then by the end of the proof we will have identified a quantity $\alpha_*>0$ to which $\|\hat\w\|$ converges; If $\alpha_*$ turns out to be independent of the choice of $K_\w$, then we may  explicitly choose $K_\w=2\alpha_*$ (say) and $\alpha_*$ is the true limit; on the other hand, if $\alpha_*$ turns out to depend on $K_\w$, this means that we could have chosen $K_\w$ arbitrarily large in the first place, and so the true limit diverges. Thus, assuming that $\|\hat\w\|$ the minimization in \eqref{eq4} is not affected by imposing the constraint $\|\w\|\leq K_\w$. With these, we can write the corresponding (AO) problem as
\begin{align}
\wt = \arg\min_{\|\w\|\leq K_\w}\max_{\|\ub\|\leq 1} \frac{1}{\sqrt{n}}\{ \|\Pp\w\|\g^T\ub - \|\ub\|\h^T\Pp\w +\ub^T(\z_\e - (\x_0^T\w)\e)	\}+\frac{\la}{\sqrt{n}}	\f(\mu\x_0+\w).\label{eq:AO1}
\end{align}
We will see that analyzing this problem is simpler than the (PO) (and certainly so of the one we started with in \eqref{eq2}). 

\subsubsection{Analysis of the Auxiliary Optimization}
The goal of this section is analyzing the (AO) problem in \eqref{eq:AO1}. In particular, we will prove i) the optimal cost of the (AO) problem converges to the optimal cost of the deterministic optimization in \eqref{eq:maxmin}, which involves three scalar optimization variables $\alpha,\beta,\tau$, ii) the max-min problem in \eqref{eq:maxmin} is \emph{strongly} convex in $\alpha$ and jointly concave in $\beta,\tau$, iii) $\|\wt\|$ converges to the unique optima $\alpha_*$ in \eqref{eq:maxmin}. With these, the claim of the Theorem follows by \cite[Thm.~1]{COLT15} (also, see \cite[Cor.~A.1]{tight}), as previously discussed.

The analysis requires several steps. The randomness in \eqref{eq:AO1} is over $\e$, $\g$, $\h$, $\x_0$ and possibly the link function $g$; at each step we condition on all but a subset of these and identify convergence of the objective function of the (AO) with respect to the remaining. Pointwise convergence (with respect to the involved optimization variables) needs to be turned into uniform convergence to guarantee that not only the objective function, but also the min/max value and the optimizer converge appropriately. (Strong) convexity of the objective will turn out to be crucial for this.

\textbf{Introducing the Frenchel conjugate}.~ To begin with, let us rewrite the (AO) problem above by expressing $f$ in terms of its Frenchel conjugate, i.e. 
\begin{align}\label{eq:fc1}
f(\x) = \sup_{\bar\vb} \bar\vb^T\x - f^*(\bar\vb)=\sup_{\bar\vb} \sqrt{n}\bar\vb^T\x - f^*(\sqrt{n}\bar\vb).
\end{align}
Translating to our problem and after rescaling this gives,
\begin{align}\label{eq:fc2}
n^{-1/2}f(\mu\x_0+\w) = \sup_{\vb} \vb^T(\mu\x_0+\w) - n^{-1/2}f^*(\sqrt{n}\vb).
\end{align}
Now, from standard optimality conditions of \eqref{eq:fc1}, the optimal $\bar\vb_*$ satisfies $\bar\vb_*\in\partial f(\x)$. Then, using condition (b) of Section \ref{sec:model}, $\|\vb_*\|=\order{\sqrt{n}}$ for all $\x$ such that $\|\x\|=\order{1}$. From this, and $\|\w+\mu\x_0\|=\order{1}$ we conclude that the optimal $\vb_*$ in \eqref{eq:fc2} satisfies $\|\vb_*\|\leq K_\vb<0$ for sufficiently large constant $K_\vb$ independent of $n$. Also, Putting everything together, \eqref{eq:AO1} is equivalent to
\begin{align}
 \min_{\|\w\|\leq K_\w}\max_{\substack{\|\ub\|\leq 1\\ 0\leq\|\vb\|\leq K_\vb }}  \frac{1}{\sqrt{n}} \ub^T(\z_\e - (\x_0^T\w)\e-
 &\|\Pp\w\|\g) -  \|\ub\|\hb^T\Pp\w \nn \\
 &+\la \vb^T(\mu\x_0+\w) - \la\fb(\vb),\label{eq:AO1.5}
\end{align}
where we have also denoted $\hb:=n^{-1/2}\h$ and $\fb(\vb)={n}^{-1/2} f^*(\sqrt{n}\vb)$. Observe again that by condition (b) of Section \ref{sec:model}, $\fb(\vb)=\max_{\x}\x^T\vb-n^{-1/2}f(\x)=\order{1}$ since $\vb=\order{1}$.

In order to somewhat simplify the exposition, we often omit explicitly carrying over the constraints $\|\w\|\leq K_\w$, $\|\vb\|\leq K_\vb$ until the very last step, but we often recall and actually make use of it.

\textbf{Optimizing over the direction of $\ub$}.~Observe that maximization over the direction of $\ub$ is easy in \eqref{eq:AO1.5}, which then becomes:
\begin{align}
 \min_{\w}\max_{\substack{0\leq\beta\leq 1\\ \vb}}  \frac{1}{\sqrt{n}} \beta\|~\z_\e - (\x_0^T\w)\e-
 &\|\Pp\w\|\g~\| -  \hb^T\Pp\w +\la \vb^T(\mu\x_0+\w) - {\la}\fb(\vb).\label{eq:AO2}
\end{align}
Now, observe that the objective function above is convex in $\w$ and jointly concave in $\beta,\vb$ (recall $f^*$ is convex). Furthermore, the constraint sets are convex and compact. Hence, we can flip the order of min-max as in \Roc:
\begin{align}
&\max_{\substack{0\leq\beta\leq 1\\\vb}}\min_{\w} \frac{\beta}{\sqrt{n}} \|~ \z_\e + (\x_0^T\w)\e - \|\Pp\w\| \g~\|  -   \hb^T\Pp\w 	+\la \vb^T(\mu\x_0+\w) - {\la}\fb(\vb) \nn \\ 
&=\max_{\substack{0\leq\beta\leq 1 \\ \vb }}\min_{\alpha_1,\alpha_2\geq 0} \frac{\beta}{\sqrt{n}} \|~ \z_\e + \alpha_2\e - \alpha_1 \g~\|  - \max_{\substack{\|\Pp\w\|=\alpha_1 \\ \x_0^T\w=\alpha_2}} \left\{  \beta \hb^T\Pp\w - \la \vb^T(\mu\x_0+\w) + \la\fb(\vb) \right\}. 
\nn
\end{align}
By decomposing $\w$ as $\Pp\w+(\x_0^T\w)\x_0$, it is not hard to perform the maximization over $\w$ to equivalently write the last display above as:
\begin{align}
\max_{\substack{0\leq\beta\leq 1 \\ \vb }}\min_{\alpha_1,\alpha_2\geq 0} \frac{\beta}{\sqrt{n}} \|~ \z_\e + \alpha_2\e - \alpha_1 \g~\|  - \alpha_1\|\beta \Pp\hb-\la\Pp\vb\| + \la\mu\vb^T\x_0 + \alpha_2\la(\vb^T\x_0) - \la \fb(\vb). 
\label{eq:AO3}
\end{align}

\textbf{The randomness of $\e$, $\g$ and $g$}.~ Until further notice condition on $\hb$ and $\x_0$. All randomness in \eqref{eq:AO3} is now on the first term. 

Consider $\beta,\vb$ fixed for now. For any pair $\alpha_1,\alpha_2$ by the WLLN, $m^{-1}\|\z_\e+\alpha_2\e-\alpha_1\g\|^2\rP\E[(g(\gamma)-\mu\gamma + \alpha_2\gamma - \alpha_1 \gamma')^2]$, where $\gamma,\gamma'\sim\Nn(0,1)$ and independent. Recall, $\E[(g(\gamma)-\mu\gamma)^2]=\sigma^2$, $\E[(g(\gamma)-\mu\gamma)\gamma]=\mu-\mu=0$ and $m/n=\delta$, to conclude that  ${n}^{-1/2}\|\z_\e+\alpha_2\e-\alpha_1\g\|\rP\sqrt{\delta}\sqrt{\sigg+\alpha_1^2+\alpha_2^2}$, where convergence is point-wise in $\alpha_1,\alpha_2$. The objective function in \eqref{eq:AO3} is jointly convex in $[\alpha_1,\alpha_2]$. Lastly, the function $\sqrt{\sigg+\alpha_1^2+\alpha_2^2}$ can be shown (by direct differentiation) to be jointly strongly convex over $[\alpha_1,\alpha_2]$. With these, we use \consist  to conclude that (for any $\beta,\vb$) i) the minimum over $\alpha_1,\alpha_2$ in \eqref{eq:AO3} converges to 
\begin{align}
\min_{\alpha_1,\alpha_2\geq 0} {\beta}\sqrt{\delta}\sqrt{\sigg+\alpha_1^2+\alpha_2^2}  - \alpha_1\|{\beta} \Pp\hb-\la\Pp\vb\| + \la\mu\vb^T\x_0 + \alpha_2\la(\vb^T\x_0) -\la \fb(\vb),
\label{eq:AO_conv0}
\end{align}
and, ii) the optimal $\alpha_1,\alpha_2$ of \eqref{eq:AO3} converge to the unique (by strong convexity) optimal of \eqref{eq:AO_conv0}. 

Up to now, $\beta,\vb$ were assumed fixed and the convergence from \eqref{eq:AO3} to \eqref{eq:AO_conv0} holds point-wise with respect to $\beta,\vb$. The objective function in \eqref{eq:AO3} is jointly concave with respect to $\beta,\vb$. Thus,
 \eqref{eq:AO3} converges to 
\begin{align}
\max_{\substack{0\leq\beta\leq 1\\ \vb}}\min_{\alpha_1,\alpha_2\geq 0} {\beta}\sqrt{\delta}\sqrt{\sigg+\alpha_1^2+\alpha_2^2}  - \alpha_1\|{\beta} \Pp\hb-\la\Pp\vb\| + \la\mu\vb^T\x_0 + \alpha_2\la(\vb^T\x_0) -\la \fb(\vb),
\label{eq:AO_conv1}
\end{align}
and the optimal $\alpha_1,\alpha_2$ of the former converge to the corresponding optima of the latter.

\textbf{Merging $\alpha_1$ and $\alpha_2$.}
It is important to note that $\alpha_1^2+\alpha_2^2$ in \eqref{eq:AO_conv1} correspond exactly to the squared norm of the error. Here, we simplify \eqref{eq:AO_conv1} by introducing the quantity $\alpha_1^2+\alpha_2^2$ as the minimization variable rather than sperately $\alpha_1$ and $\alpha_2$. By first order optimality conditions in \eqref{eq:AO_conv1} we find
\begin{align}\label{eq:alphas}
\alpha_1 \beta\sqrt{\delta} =  \|\beta\Pp\hb-\la\Pp\vb\|\sqrt{\alpha_1^2+\alpha_2^2+\sigg}\quad\text{ and }\quad \alpha_2 -\beta\sqrt{\delta} =\la\vb^T\x_0 \sqrt{\alpha_1^2+\alpha_2^2+\sigg} .
\end{align}
Substituting this in \eqref{eq:AO_conv1}, the objective becomes (ignoring the terms that do not involve $\alpha_1$ or $\alpha_2$):
\begin{align}\nn
 {\beta}\sqrt{\delta}\sqrt{\sigg+\alpha_1^2+\alpha_2^2}  - \frac{\sqrt{\sigg+\alpha_1^2+\alpha_2^2}}{\beta\sqrt{\delta}}\left( \|\beta\Pp\hb-\la\Pp\vb\|^2 + (\la\vb^T\x_0)^2  \right)
\end{align}
But, from \eqref{eq:alphas} we find $\sqrt{\sigg+\alpha_1^2+\alpha_2^2}\sqrt{\|\beta\Pp\hb-\la\Pp\vb\|^2 + (\la\vb^T\x_0)^2} = \beta\sqrt{\delta}\sqrt{\alpha_1^2+\alpha_2^2}$. Combining, we conclude that \eqref{eq:AO_conv1} can be written as
\begin{align}
\max_{\substack{0\leq\beta\leq 1\\ \vb}}\min_{\alpha\geq 0} {\beta}\sqrt{\delta}\sqrt{\sigg+\alpha^2}  - \alpha\|\beta\Pp\hb-\la\vb\|  + \la\mu\vb^T\x_0 -\la \fb(\vb),
\label{eq:AO_conv2}
\end{align}
where the new optimization variable $\alpha$ plays the role of $\sqrt{\alpha_1^2+\alpha_2^2}$, thus it represents the norm of the error vector $\|\w\|$. We have also identified $\|\beta\Pp\hb-\la\vb\|^2 + (\la\vb^T\x_0)^2 = \|\beta\Pp\hb-\la\vb\|^2$


\textbf{Introducing a new optimization variable}. To get a better handle at it, we square the norm term in \eqref{eq:AO_conv2} at the expense of introducing a new scalar optimization variable. This is based on the following trick: 
\begin{align}\label{eq:trick}
\sqrt{x}=\min_{\tau>0}\frac{\tau}{2}+\frac{x}{2\tau}
\end{align}
 for any $x\geq 0$. Thus, \eqref{eq:AO_conv2} becomes
\begin{align}
\max_{\substack{0\leq\beta\leq 1\\ \vb, \tau> 0}}\min_{\alpha\geq 0} {\beta}\sqrt{\delta}\sqrt{\sigg+\alpha^2}  - \frac{\alpha \tau}{2} - \frac{\alpha}{2\tau}{\|\beta\Pp\hb-\la\vb\|^2}  + \la\mu\vb^T\x_0 -\la \fb(\vb),
\label{eq:AO_conv3}
\end{align}
where we have also flipped the order of min-max between $\alpha$ and $\tau$. We could do this as in \Roc since the objective is convex in $\alpha$ and concave in $\tau$, the constraint sets are both convex and both of them are bounded. To argue the boundedness, recall that $\alpha\leq K_\w$; for $\tau$ it suffices to combine optimality conditions of \eqref{eq:trick} and boundedness of $\vb$, $\|\vb\|_2\leq K_\vb$.

\textbf{Optimizing over $\vb$.} Note that the objective in \eqref{eq:AO_conv3} is concave in $\vb$, convex in $\alpha$ and the constraint sets are convex compact. Thus, as it might be expected by now, we use \Roc to flip the corresponding order of max-min. Also, after some simple algebra while using $\Pp\x_0=0$ and $\|\x_0\|=1$, it can be shown that 
\begin{align}\nn
{\|\beta\Pp\hb-\la\vb\|^2 }  - 2\frac{\tau}{\alpha}\la\mu\vb^T\x_0 = \|\la\vb - (\beta\Pp\hb+\frac{ \tau}{\alpha}\mu\x_0)\|^2 - \mu^2\frac{\tau^2}{\alpha^2}.
\end{align}
Combining, we conclude with 
\begin{align}
\eqref{eq:AO_conv3}&=\max_{\substack{0\leq\beta\leq 1\\  \tau> 0}}\min_{\alpha\geq 0} {\beta}\sqrt{\delta}\sqrt{\sigg+\alpha^2}  - \frac{\alpha p}{2} + \mu^2\frac{\tau}{2\alpha} - \frac{\alpha\la^2}{\tau}\min_{\vb}\left\{ \frac{1}{2}\|\vb - (\frac{\beta}{\la}\Pp\hb+\frac{ \tau}{\alpha\la}\mu\x_0)\|^2 + \frac{\tau}{\la\alpha} \fb(\vb) \right\}.
\label{eq:AO_conv4} \\&= \max_{\substack{0\leq\beta\leq 1\\  \tau> 0}}\min_{\alpha\geq 0} G(\alpha,\beta,\tau) \nn
\end{align}

Here, ${G}(\al,\beta,\tau)$ is convex in $\al$ (see \eqref{eq:AO_conv2}) and jointly concave in  over $\beta$, $\tau$. To see the latter it suffices to show that $ \frac{\al\la^2}{\tau}\| \vb	-	\frac{\beta}{\la}( \Pp\hb+\frac{\mu\tau}{\la\al}\x_0)\|^2$ is jointly convex over $\beta,\tau,\vb$ (minimization over $\vb$ does not change the joint convexity over $\tau$ and $\beta$.). Norm is separable over its entries, so we equivalently show that for scalars $\tau,\beta,v$, the function $\frac{1}{\tau}(v-c_1\beta-c_2\tau)^2$ is jointly convex over $\tau>0,\beta$; this is true as the perspective function of $(v-c_1\beta-c_2)^2$. 

\textbf{The randomness of $\hb$ and $\x_0$}.~ For now, fix any $\beta,\tau$ and let $\ah:=\ah(\beta,\tau)$ be the minimizer in \eqref{eq:AO_conv4}.

First, we prove that $\ah(\beta,\tau)>0$. For any $\alpha\geq 0$, by choosing $\tilde\vb=\min\{{\frac{\tau\mu}{\alpha\la},\frac{K'_\vb}{\|\x_0\|}}\}\x_0=:\theta\x_0$ where $0<K'_\vb\leq K_\vb$ such that $\tilde\vb\in \dom\fb$, we find
\begin{align}
&\min_{\vb}\left\{ \frac{{\alpha\la^2}}{2\tau}\|\vb - (\frac{\beta}{\la}\Pp\hb+\frac{ \tau}{\alpha\la}\mu\x_0)\|^2 + \la \fb(\vb) \right\} \leq \frac{{\alpha\la^2}}{2\tau}\|\tilde\vb - (\frac{\beta}{\la}\Pp\hb+\frac{ \tau}{\alpha\la}\mu\x_0)\|^2 + \la \fb(\tilde\vb) \nn \\
&\qquad\qquad\qquad\quad\leq \frac{\mu^2\tau}{2\alpha}\left( 1-\frac{\theta\alpha\la}{\tau\mu} \right)^2 + \frac{\alpha}{2\tau}\beta^2\|\Pp\h\|^2+\la\fb(\theta\x_0).
\end{align}
Thus, the value of the objective in \eqref{eq:AO_conv4} is lower bounded by
$$
\beta\sqrt{\delta}\sqrt{\sigg+\alpha^2}  - \frac{\alpha p}{2} + \mu^2\frac{\tau}{2\alpha} \left(1 - \left( 1-\frac{\theta\alpha\la}{\tau\mu} \right)^2 \right) - \frac{\alpha}{2\tau}\beta^2\|\Pp\h\|^2- \la\fb(\theta\x_0).
$$
which goes to $+\infty$ as $\alpha\rightarrow 0$, since by definition $\theta\leq\frac{\tau\mu}{\alpha\la}$. Hence, $\ah>0$, as desired.

Now, fix $\beta,\tau,\alpha>0$, denote $c_1=\frac{\beta}{\la}, c_2=\frac{\tau}{\alpha\la}, c_3=\frac{\tau\mu}{\alpha\la}$ and consider
$$
R(\hb,\x_0):=R(\alpha,\beta,p;\hb,\x_0):= \min_{\vb}\left\{ \frac{1}{2}\|\vb - c_1\Pp\hb-c_2\x_0\|^2 + c_3 \fb(\vb) \right\}.
$$
 Recall from Assumption 1 that  
\begin{align}\label{eq:A}
A(\hb,\x_0):= R(\alpha,\beta,p;\hb,\x_0):= \min_{\vb}\left\{ \frac{1}{2}\|\vb - c_1\hb-c_2\mu\frac{\xx_0}{\sqrt{n}}\|^2 + c_3 \fb(\vb) \right\} 
\end{align}
converges to $F:=F(c_1,c_2,c_3)$ in probability. Also, recall $\xx_0=\x_0\|\xx_0\|$. Next, we show that for all constant $\zeta>0$
\begin{align}\label{eq:2prove}
| R(\hb,\x_0) - A(\hb,\x_0) | \leq \zeta
\end{align}
with probability approaching one in the limit of $n\rightarrow\infty$. Combining this with Assumption 1, will prove that $A(\hb,\x_0)$ converges in $F$ in probability.

\underline{Proof of \eqref{eq:2prove}}: Fix any $\eps>0$. We condition on the following events: 
\begin{align}\label{eq:events}
\begin{cases}
~|\hb^T\x_0|\leq \eps, \\
~1-\eps\leq {n^{-1/2}\|\xx_0\|}\leq 1+\eps.
\end{cases}
\end{align}
Each one of the events occurs with probability approaching one as $n\rightarrow\infty$; the first follows since $\hb\sim\Nn(0,n^{-1/2})$ and $\|\x_0\|=1$ and from standard tail bounds on Gaussians; the second is due to condition (b) of Section \ref{sec:model}.
Without loss of generality assume $R(\hb,\x_0)\geq A(\hb,\x_0)$, and let $\vbs$ be optimal in \eqref{eq:A}, then
\begin{align}
&| R(\hb,\x_0) - A(\hb,\x_0) | \leq  \frac{1}{2}\|\vbs - c_1\Pp\hb-c_2\x_0\|^2 - \frac{1}{2}\|\vbs - c_1\hb-c_2\frac{\xx_0}{\sqrt{n}}\|^2 \nn \\
&= \left(c_1(\x_0^T\hb)\x_0+c_2\xx_0\left(\frac{1}{\sqrt{n}}-\frac{1}{\|\xx_0\|}\right)\right)^T
\left( \vb_* - c_1\hb - \frac{1}{2} c_2\xx_0\left(\frac{1}{\sqrt{n}}+\frac{1}{\|\xx_0\|}\right) + \frac{1}{2}c_1(\x_0^T\hb)\x_0 \right) \nn\\
&= - \frac{1}{2}c_1^2(\x_0^T\hb)^2 + c_1(\x_0^T\h)(\x_0^T\vbs)  - c_1c_2(\x_0^T\hb) \frac{\|\xx_0\|}{\sqrt{n}}
 + c_2(\x_0^T\vbs)\left(\frac{\|\xx_0\|}{\sqrt{n}}-1\right) -  \frac{1}{2}c_2^2\left(\frac{\|\xx_0\|^2}{{n}}-{1}\right)\nn\\
& \leq \frac{1}{2}c_1^2\eps^2 + c_1 \|\vbs\|\eps + c_1c_2\eps(1+\eps) + c_2\|\vbs\|\eps + \frac{1}{2}c_2^2\eps(2+\eps)
\end{align}
where the last line follows after bounding the absolute values of the summands using \eqref{eq:events}. Recall now that $\|\vbs\|\leq K_\vb<\infty$ and also $c_1,c_2,c_3$ are also bounded constants (independent of $n$). Then, for all $\zeta>0$ in \eqref{eq:2prove} we can find sufficiently small $\eps>0$ such that the value of the last expression in the panel above is no larger than $\zeta$, thus completing the proof of \eqref{eq:2prove}

Thus, we have shown that $G(\alpha,\beta,\tau)$ in \eqref{eq:AO_conv4} converges pointwise to
$$H(\alpha,\beta,\tau) := \sqrt{\delta}\sqrt{\sigg+\alpha^2}  - \frac{\alpha p}{2} + \mu^2\frac{\tau}{2\alpha} - \frac{\alpha\la^2}{\tau} F(\frac{\beta}{\la},\frac{\tau}{\alpha\la},\frac{\tau\mu}{\alpha\la}),$$
in the limit of $n\rightarrow\infty$. Note that $H$ is strongly convex in $\alpha$ and jointly concave in $\beta,\vb$ since taking limits does not affect convexity properties (recall that $G$ is convex-concave). Also, we showed that  $\alpha_*(\h,\x_0)>0$ for the optimal in \eqref{eq:AO_conv4}. With these, it follows as per \consist that (i) 
\begin{align}\label{eq:almost?}
\min_{0\leq\beta\leq 1,\tau>0}\max_{\alpha>0}G(\alpha,\beta,\tau)\rP\min_{0\leq\beta\leq 1,\tau>0}\max_{\alpha>0}H(\alpha,\beta,\tau),
\end{align}
and, (ii) $\alpha_*(\h,\x_0)\rP\alpha_*$, where $\alpha_*$ the unique minimizer of the second optimization in \eqref{eq:almost?}. This completes the proof of the Theorem. 

\subsection{Theorem \ref{coro:main}}
The theorem is a direct consequence of Theorem \ref{thm:main}. In particular, Theorem \ref{thm:main} proves that the value $\alpha_*$ to which the error converges only depends on $g$ through the parameters $\mu$ and $\sigg$. Those are the same (by definition) for the non-linear and the linear case considered, thus the errors are the same.

\appendix[Proof of Theorem \ref{coro:sparse}]\label{sec:sparse_app}

Specializing Theorem \ref{thm:main} to the setup of Section \ref{sec:sparse} we showed in the same section that $\|\hat\x-\mu\x_0\|$ converges in probability to the unique minimizer $\alpha_*$ of the following max-min problem:
\begin{align}\label{eq12}
\max_{\substack{0\leq\beta\leq 1\\ \tau> 0}}\min_{\al\geq 0}  H(\al,\beta,\tau):=
		{\beta}\sqrt{\delta}\sqrt{\sigma^2+\al^2} 	-\frac{\tau\alpha}{2}	+\frac{\tau\mu^2}{2\al}- \frac{\al}{2\tau}	\E\left[\eta^2\left({\beta} h+\frac{\mu\tau}{\al}\XX_0;\la\right)\right],
\end{align}
where the expectation is over $h\sim\Nn(0,1)$ and $\XX_0\sim p_{\XX_0}$. Here, we prove Theorem \ref{thm:main} by analyzing the optimality conditions of \eqref{eq12}. Recall as in the IEEEproof of Theorem \ref{thm:main} that $H$ is jointly concave in $\beta,p$ and strongly convex in $\alpha$. 

\subsection{First Order Optimality Conditions.}
We begin with a lemma, which characterizes the first-order optimality conditions of \eqref{eq12}.
\begin{lem}[Optimality Conditions]\label{lem:foc}
Consider the following pair of equations with respect to $\beta$ and $\kappa$:
\begin{numcases}{}
\beta^2\kappa^2\delta = \sigma^2 + \E\left[( \eta(\beta\kappa h + \mu \XX_0; \kappa\la) - \mu\XX_0 )^2\right],\label{eq:set1}\\
\beta\kappa\delta = \E[( \eta(\beta\kappa h + \mu \XX_0; \kappa\la)\cdot h )].\label{eq:set2}
\end{numcases}
Also, define $\lam$ to be the unique non-negative solution to the equation
$$
(1+x^2)\int_{-\infty}^{-x}e^{-z^2/2}\mathrm{d}z-xe^{-x^2/2}=\delta\sqrt{\frac{\pi}{2}}.
$$
With these, let $(\beta_*,\tau_*,\alpha_*)$ be optimal in \eqref{eq12}. Then,
\begin{align}\label{eq:alka}
\alpha_*^2 = \beta_*^2\kappa_*^2\delta - \sigg
\quad \text{ and } \quad \kappa_*=\frac{\sigma}{\sqrt{\beta_*^2\delta-\tau_*^2}}.
\end{align}
such that,
\begin{enumerate}[(i)]
\item If $\beta_*=1$ and $\la>\lam$, then $\kappa_*$ is the unique solution to \eqref{eq:set1} for $\beta=1$,
\item If $\beta_*\in(0,1)$ ,then $\kappa_*,\beta_*$ are solutions to the pair of equation \eqref{eq:set1}-\eqref{eq:set2}.
\end{enumerate}
\end{lem}

\begin{IEEEproof}
Let us compute $\frac{\partial}{\partial \al}H(\beta,\al,\tau)$ and $\frac{\partial}{\partial \tau}H(\beta,\al,\tau)$. For convenience  define 
$$P\left(\frac{\tau}{\al}\right):=\frac{\tau\mu^2}{2\al}-\frac{\al\la^2}{2\tau}	\E\left[\eta^2\left(\frac{\beta}{\la} h	+\frac{\mu\tau}{\la\al}\XX_0;1\right)\right].
$$ 
Taking derivatives in \eqref{eq12} with respect to $\alpha$ and $\tau$ and equating them with zero gives
\begin{subequations}\label{Der}
\begin{align}
\frac{\beta\sqrt{\delta}\alpha}		{\sqrt{\al^2+\sigma^2}}		-\frac{\tau}{2}		-\frac{\tau}{\al^2}P'(\frac{\tau}{\al})&=0,\label{Der1}\\
-\frac{\al}{2}+\frac{1}{\al}P'(\frac{\tau}{\al})&=0.\label{Der2}
\end{align}
\end{subequations}
Here, $P'$ is the derivative of $P(x)$ with respect to $x$. Any optimal $\beta_*,\tau_*,\alpha_*$ satisfies these. Then, it only takes multiplying \eqref{Der2} by $\frac{\tau}{\al}$ and adding the result to \eqref{Der1} to see that
\begin{align}\label{eq:rel}
\alpha_*=\frac{\tau_*\sig}{\sqrt{\beta^2\delta-\tau_*^2}}.
\end{align}

Next, substituting \eqref{eq:rel} in \eqref{Der2} it can be shown that,
\begin{align}
-\frac{\sigma^2}{2}+\frac{\sig^2 }{2\tau^2}	\E[\big(\eta(\beta h +\frac{\sqrt{\beta^2\delta-\tau^2}}{\sig}\mu\XX_0;\la)-\frac{\sqrt{\beta^2\delta-\tau^2}}{\sig}\mu\XX_0\big)^2]=0\nn.
\end{align}

To reach this we have also used the following facts: $\eta(x;\la)\frac{\partial}{\partial x}\eta(x;\la)=\eta(x;\la)$, $ \la\eta(\frac{x}{\la};1)=\eta(x;\la)$ and $\E[\XX_0^2]=1$ by assumption. Multiplying the result with $ {2\tau^2}/{\sig^2} $ and defining $$\kappa:=\frac{\sigma}{\sqrt{\beta^2\delta-\tau^2}},$$ we conclude with,
\begin{align}
\beta^2\delta\kappa^2-\sig^2=\E[\big( \eta(\beta\kappa h +\mu \XX_0;\kappa\la)-\mu\XX_0 \big)^2]\label{eq14},
\end{align}
which is same as \eqref{eq:set1}.
Also, with respect to the optimal $\kappa_*$ it is easily seen by \eqref{eq:rel} that
\begin{align}\label{eq:alphakappa}
\alpha_*^2 = \beta_*^2\kappa_*^2\delta - \sigg.
\end{align}
%
The derivative in \eqref{eq12} with respect to $\beta$ gives
\begin{align}
\frac{\partial}{\partial \beta}H(\al,\beta,\tau)&=\sqrt{\delta}\sqrt{\sig^2+\al^2}-\frac{\alpha}{\tau}\E[\eta(\beta h+\frac{\mu\tau}{\alpha}\XX_0;\la)h]\nonumber\\
=&\beta\delta\kappa-\kappa\E[\eta(\beta h+\frac{\mu\XX_0}{\kappa},\la)h]=\beta\delta\kappa-\E[\eta(\kappa\beta h+\mu\XX_0;\la\kappa)h].\label{eq:bder}
\end{align}
where we have also used  \eqref{eq:alphakappa}. Note that the above is same as \eqref{eq:set2} and recall the constraint $0\leq\beta\leq 1$ in \eqref{eq12} to conclude with the desired.

It only remains to show that the solution with respect to $\kappa$ of \eqref{eq:set1} (eqv. of \eqref{eq14}) is unique when $\beta=1$ and $\la\geq\lam$. For $\beta=1$, \eqref{eq:set1} is the same as fixed point equation \cite[Eqn.~(1.9)]{montanariLASSO}, which in turn was shown to admit a unique solution for all $\la>\lam$ in \cite{DMM} (see \cite[Prop.~1.3]{montanariLASSO}).
\end{IEEEproof}

\subsection{The Regions of Operation}

We build up to the IEEEproof of Theorem \ref{coro:sparse} through a series of auxiliary lemmas. Through the lemmas, we identify two ``regimes of operation" of the LASSO. The first, we call $\Rbad$, and it corresponds to values of $\la$ for which the optimal $\beta$ is in the open set $(0,1)$. The second regime, is such that $\beta=1$. If $\delta<1$, we prove in Lemma \ref{lem:<1} that there exists a unique critical value $\lac$ separating the two regimes in the sense that $\Rbad$ extends from $0$ to $\lac$. If on the other hand $\delta\geq 1$, then there is no $\Rbad$ region  (Lemma \ref{lem:>1}).

First, we need a few useful definitions.

\begin{defn}
For any $\la>0$, we let $\alpha_*(\la)$, $\tau_*(\la)$ and $\beta_*(\la)$ be  optimal solutions in \eqref{eq12}. Apart from $\alpha_*(\la)$, the others are not necessarily unique at this point. Also, $\kappa_*(\la)$ is defined as in \eqref{eq:alka}.
\end{defn}

\begin{defn}[Bad Regime]
We say that a value $\la>0$ is in the bad regime $\Rbad$, denote $\la\in\Rbad$, if there exists $\beta_*(\la)\in(0,1)$.
\end{defn}

\begin{defn}[Critical Regime] We say that  a value $\lac>0$ is in the critical regime $\Rcrit$, denote $\lac\in\Rcrit$ if for some $\kac$, the pair $\lac,\kac$ solves:
\begin{numcases}{}
\kappa^2\delta = \sigma^2 + \E\left[( \eta(\kappa h + \mu \XX_0; \kappa\la) - \mu\XX_0 )^2\right],\label{eq:tau_app}\\
\kappa\delta = \E[( \eta(\kappa h + \mu \XX_0; \kappa\la)\cdot h )].\label{eq:beta_app}
\end{numcases}
\end{defn}

As an immediate consequence of the definition above and the first order optimality conditions in Lemma \ref{lem:foc}, we have 
\begin{align}\label{eq:lac_use}
\beta_*(\lac)=1, \quad \kappa_*(\lac)=\kac \quad \text{and} \quad \alpha_*(\lac) = \sqrt{\delta\kac^2-\sigg}.
\end{align}

Also, the following lemma reveals the importance of $\lac$: all $\la<\lac$ are in $\Rbad$ and the squared error is constant in that regime, i.e. $\alpha_*(\la)=\alpha_*(\lac)$.

\begin{lem}[Error in $\Rbad$]\label{lem:Rbad}
Let $\lac\in\Rcrit$. Then, for all $0<\la'<\lac$, it holds $\la'\in\Rbad$. Furthermore, 
$\beta_*(\la')=\la/\lac$, $\la'\kappa_*(\la')=\kac\lac$ and $\alpha_*(\la')=\alpha_*(\lac)$.
\end{lem}
\begin{IEEEproof}
Fix any $0<\la'<\lac$. By definition, there exists $\kac$ such that $\lac,\kac$ satisfy \eqref{eq:tau_app}-\eqref{eq:beta_app}. Define $\beta':=\la/\lac$ and $\kappa':=\kac/\beta'$. It is then easy to see that $\beta',\kappa'$ solve \eqref{eq:set1}-\eqref{eq:set2} (for $\la=\la'$ therein). Also, $\beta'<1$ by definition. Thus, $\la'\in\Rbad$ and $\beta_*(\la')=\la/\lac, \kappa_*(\la')=\kac\lac/\la'$. Also, using \eqref{eq:alka} and \eqref{eq:lac_use}, $\alpha_*(\la) = \sqrt{\delta\beta^2_*(\la')\kappa^2_*(\la)-\sigg}=\sqrt{\delta\kappa_*^2(\lac)-\sigg}=\alpha_*(\lac)$.
\end{IEEEproof}

It is thus important to identify the critical values of the regularizer parameter, i.e. all $\lac\in\Rcrit$. Values in $\Rbad$ are important towards this direction, since as shown in the next lemma, for any $\la\in\Rbad$ there must exist some $\lac>\la$.

\begin{lem}[$\Rbad\rightarrow\lac$]\label{lem:right}
Let $\la_1\in\Rbad$, then there exists $\la_2\in\Rcrit$ with $\la_2>\la_1$. 
\end{lem}
\begin{IEEEproof}
Let $\beta_1,\alpha_1,\kappa_1$ be optimal corresponding to $\la_1$. Since $\la_1\in\Rbad$, it holds $0<\beta_1< 1$. Then, from Lemma \ref{lem:foc}, $\kappa_1,\beta_1$ solve \eqref{eq:set1}-\eqref{eq:set2}. Starting from these and substituting $\la_2:=\la_1/\beta_1$ and $\kappa_2:=\kappa_1\beta_1$ therein, it is not hard to see that this is equivalent with $\la_2,\kappa_2$ satisfying \eqref{eq:tau_app}-\eqref{eq:beta_app}. Thus, $\la_2\in\Rcrit$. Also, clearly $\la_2>\la_1$. 
\end{IEEEproof}

The lemma below is important since it shows that when $\delta<1$ there exists a \emph{unique} $\lac\in\Rcrit$.
\begin{lem}[Unique $\lac$]\label{lem:<2}
Suppose $\delta<1$. The set of equations \eqref{eq:tau_app}-\eqref{eq:beta_app} has a \emph{unique} pair of solutions $(\kappa,\la)$. Thus, there exists unique $\lac\in\Rcrit$.
\end{lem}
\begin{IEEEproof} 

First, we show that there exists at most one $\lac\in\Rcrit$. For the shake of contradiction assume two different pairs of solutions, say $(\kappa_1,\la_1)$ and $(\kappa_2,\la_2)$. By definition, $\la_1,\la_2\in\Rcrit$. First, note that we cannot have $\la_1=\la_2$, since if this was the case then from \eqref{eq:lac_use} we would also have $\kappa_1=\kappa_2$. Henceforth, assume w.l.o.g. that $\la_1<\la_2$. It follows from Lemma \ref{lem:Rbad} that $\la_1\in\Rbad$ and also $\kappa_*(\la_1)\la_1 = \kappa_*(\la_2)\la_2$. Thus, 
\beq\label{eq:2contra1}
\kappa_*(\la_1)<\kappa_*(\la_2).
\eeq
But also, again from Lemma \ref{lem:Rbad}, $\alpha_*(\la_1)=\alpha_*(\la_2)$. Since, $\la_1,\la_2\in\Rcrit$, this implies when combined with \eqref{eq:lac_use} that $\kappa_*(\la_1)=\kappa_*(\la_2)$, which contradicts \eqref{eq:2contra1}, completing the IEEEproof of this part.

Let us now prove that $\Rcrit$ is non-empty.
To begin with, we show that $\Rbad$ is non-empty in this case. In particular, we show that $\lam$ defined in Lemma \ref{lem:foc} is in $\Rbad$. Since, $\delta<1$, we have $\lam>0$. Suppose that ($\beta_*(\lam)=1,\kappa_*(\lam))$ is optimal for some $\kappa_*(\lam)$, then, from first-order optimality conditions, $\kappa_*(\lam),\lam$ solves \eqref{eq:set1} for $\beta=1$. But, then as in \cite[pg.~16]{montanariLASSO} $\kappa_*(\lam)\rightarrow\infty$. Also, since $H(\alpha,\tau,\beta)$ is concave in $\beta$, the above imply that $\frac{\partial H}{\partial\beta}\big|_{(\beta=0,\kappa\rightarrow\infty)}\geq 0$, or equivalently from \eqref{eq:bder}, 
$$
\int_{\lam}^{\infty} h(h-\lam)e^{-h^2/2}\mathrm{d}h \leq {\delta}\sqrt\frac{\pi}{2}.
$$
Recalling the definition of $\lam$ in Lemma \ref{lem:foc}, it can be shown (using standard inequalities on tail functions of gaussians) that the inequality above is violated for all $0<\delta<1$. Hence, it must be $\beta_*(\lam)<1$. Also, $\beta_*(\lam)>0$ because of \eqref{eq:rel}.  Thus, $\lam\in\Rbad$. 
To complete, the IEEEproof use Lemma \ref{lem:right} with $\la_1=\lam$ to see that there exists $\la_2\in\Rcrit$.
\end{IEEEproof}

\begin{lem}[$\delta<1$]\label{lem:<1}
Suppose $\delta<1$ and let $\lac\in\Rcrit$. Furthrermore, i) for all $\la\leq\lac$, $\alpha_*(\la)=\alpha_*(\lac)$, and, ii) for all $\la>\lac$, $\kappa_*(\la)$ is the unique solution to \eqref{eq:set1} for $\beta=1$.
\end{lem}
\begin{IEEEproof}
Existence and uniqueness of $\lac$ is proved in Lemma \ref{lem:<2}

i) For $\la\leq\lac$, the claim follows directly from Lemma \ref{lem:Rbad}. 

ii) Next, we show that for $\la\geq\lac$, there exists an optimal solution for which $\beta_*(\la)=1$. This suffices since then $\kappa_*(\la)$ is indeed solving \eqref{eq:set1} for $\beta=1$ (by first order  optimality conditions), and, also, the solution is unique by \cite{DMM},\cite[Prop.~1.3]{montanariLASSO} and the fact that $\lam\leq\lac\leq\la$. 
To see that $\beta_*(\la)=1$, we argue as follows. First, $\beta_*(\la)\not\in(0,1)$. Otherwise, $\la\in\Rbad$, thus, by Lemma \ref{lem:right} there exists $\la'>\la\geq\lac$ such that $\la'\in\Rcrit$, which contradicts the uniqueness of $\lac$. Hence, $\beta_*(\la)=1$.
\end{IEEEproof}

\begin{lem}[$\delta>1$]\label{lem:>1}
Suppose $\delta> 1$, then for all $\la\geq 0$, $\kappa_*(\la)$ is the unique solution to \eqref{eq:set1} for $\beta=1$.
\end{lem}
\begin{IEEEproof}
First, let us show that for $\la\rightarrow0$, the optimal $\beta_*(\la)=1$. Indeed for $\beta=1$ and $\la\rightarrow\infty$, \eqref{eq:bder} gives
$$
\frac{\partial H}{\partial\beta} = \delta - \E[(h+\frac{\mu}{\kappa}\XX_0)h] = \delta - 1 > 1.
$$
Thus, from concavity of $H$ with respect to $\beta$, we find that the unique optimal value for $\beta$ is 
\begin{align}\label{eq:2contra}
\beta_*(\la\rightarrow 0)=1.
\end{align}
 Also, as in the IEEEproof of Lemma \ref{lem:<1}, $\beta_*(\la\rightarrow \infty)=1$. Thus, again similar to Lemma \ref{lem:<1}, it suffices to prove that there exists no $\la\in\Rbad$. For the shake of contradiction, suppose that there exists $\la_1\in\Rbad$. By Lemma \ref{lem:right}, there exists $\la_1<\lac\in\Rcrit$. But, then $\beta_*(\la\rightarrow 0)\rightarrow 0$, which contradicts \eqref{eq:2contra}. This completes the IEEEproof.
\end{IEEEproof}

\begin{IEEEproof}(of Theorem \ref{coro:sparse}) The claim of the theorem is now a direct consequence of Lemmas \ref{lem:<1} and \ref{lem:>1} combined with \eqref{eq:alphakappa}.
\end{IEEEproof}

\appendix[Proofs for Section \ref{sec:qbit}]\label{sec:qbit_app}

\subsection{The LM Algorithm}
The Lloyd-Max algorithm is an algorithm for finding the quantization threshold ${t_i}$ and the representation
points ${\ell_i}$. Given real values $x\in\R$ sampled from some probability density $\phi(x)$ it looks for optimal sets $\hat\tb$, $\hat\ellb$ that minimizes the mean-square-error (MSE) between $x$ and their corresponding quantized values $Q_q(x;\ellb,\tb)$, i.e.
\begin{align}\label{eq:MSE}
(\hat\ellb,\hat\tb) := \arg\min_{\ellb,\tb} \E_{x\sim\phi}[ ( x - Q_q(x;\ellb,\tb) )^2 ].
\end{align}
 The algorithm simply alternates between i) optimizing the threshold $t_i$ for a given set of $\ellb$, and then ii) optimizing the levels $\ellb_i$ for the new thresholds. It is well known that the converging points $\ellbLM,\tbLM$ of the algorithm satisfy
\begin{subequations}\label{eq:LM}
\begin{align}
&\tbLM_i=\frac{\ellbLM_i+\ellbLM_{i+1}}{2}&\qquad &i=1,..., L-1,\\
&\ellbLM_i=\left({\int_{\tbLM_{i-1}}^{\tbLM_i}\phi(x)\mathrm{d}x}\right)^{-1}\left({\int_{\tbLM_{i-1}}^{\tbLM_i}x\phi(x)\mathrm{d}x}\right)&\qquad &i=1,...,L.
\end{align}
\end{subequations}
Furthermore, they are stationary points of the objective function in \eqref{eq:MSE}.

\subsubsection{Gaussian case}

Assume that the values $x$ are sampled from a standard gaussian distribution, i.e. $x\sim\Nn(0,1)$ and $\phi(x)=(1/\sqrt{2\pi})\exp(-x^2/2)$. Also, recall the definition of the parameters $\mu,\sigg$ in \eqref{eq:musig2}; setting $g=Q_q$ therein, we find 
\begin{subequations}\label{eq:mut}
\begin{align}
\mu:= \mu(\ellb,\tb)=2\sum_{i=1}^L \ell_i\int_{t_{i-1}}^{t_i}x\phi(x)\mathrm{d}x\\
\tau^2:= \tau^2(\ellb,\tb)=2\sum_{i=1}^L\ell_i^2 \int_{t_{i-1}}^{t_i}\phi(x)\mathrm{d}x
\end{align}
\end{subequations}
In this notation, the objective in \eqref{eq:MSE} can be writthen as $\tau^2-2\mu+1$. Thus, $\ellbLM,\tbLM$ satisfy
\begin{align}\label{eq:LMopt}
(\tau^2)'\big|_{(\ellbLM,\tbLM)}  = 2\mu'\big|_{(\ellbLM,\tbLM)},
\end{align}
Here and onwards we use $(\tau^2)'$, $\mu'$ to denote the gradient of $\tau^2$ and $\mu$ with respect to the vector $[\ellb^T,\tb^T]$. The gradients are evaluated at the point $(\ellbLM,\tbLM)$ in \eqref{eq:LMopt}

\subsection{q-Bit Compressive Sensing}
We prove that the LM algorithm is an efficient algorithm when the objective is minimizing the LASSO reconstruction error of a signal $\x_0$ to which we have access through $q$-bit quantized linear measuments $Q_q(\ab_i^T\x;\ellb,\tb)$. It was shown in Section \ref{sec:conc} that the problem can be posed as that of finding $\ellb_*,\tb_*$ such that
\begin{align}\label{eq:Qobj_app}
(\tb_*,\ellb_*) = \arg\min_{\tb,\ellb} \frac{\sigg(\tb,\ellb)}{\mu^2(\tb,\ellb)}=\arg\min_{\tb,\ellb} \frac{\tau^2(\tb,\ellb)}{\mu^2(\tb,\ellb)}.
\end{align}
The following Lemma proves the claim made in Section \ref{sec:claim}, i.e. the converging point of the LM algorithms are stationary points of the objective function in \eqref{eq:Qobj_app}.

\begin{lem}
Then, the converging points of the LM algorithm, say $(\tbLM,\ellbLM)$ satisfy
\begin{align}
&\frac{\partial}{\partial \ell_i}\left(\frac{\tau^2(\ellb,\tb)}{\mu^2(\ellb,\tb) } \right)\bigg|_{(\ellb,\tb)=(\ellbLM,\tbLM)}=0&,\qquad& i=1,...,L,\nonumber\\
&\frac{\partial}{\partial t_i}\left(\frac{\tau^2(\ellb,\tb)}{\mu^2(\ellb,\tb) }\right)\bigg|_{(\ellb,\tb)=(\ellbLM,\tbLM)}=0,&\qquad& i=0,...,L-1.\label{deriv}
\end{align}
\end{lem}
\begin{IEEEproof}
Call $R(\tb,\ellb) = \frac{\tau^2(\tb,\ellb)}{\mu^2(\tb,\ellb)}$. We denote $R':=R'(\tb,\ellb)$ for its gradient with respect to the vector $[\tb^T,\ellb^T]$. It suffices to prove that $R'\big|_{(\ellbLM,\tbLM)}=0$, or equivalently, that at the point $(\tb,\ellb)=(\tbLM,\ellbLM)$ the following holds:
\begin{align}\label{eq:want}
(\tau^2)'\mu^2 = 2\tau^2\mu\mu'.
\end{align}
To see that this is the case,  note that
\begin{align}\label{eq:first}
\tau^2(\tbLM,\ellbLM) = \mu(\tbLM,\ellbLM)
\end{align}
This follows by direct substitution of combining  \eqref{eq:LM} in \eqref{eq:mut}.
Then, \eqref{eq:want} follows from \eqref{eq:first} and \eqref{eq:LMopt}.
\end{IEEEproof}

\end{document}